\documentclass[10 pt]{article}

\font\fiverm=cmr5 
\catcode`@=11 \catcode`!=11

\@ifundefined{normalshape}{%
     \def\pdollar{{\ifdim \fontdimen\@ne\font >\z@ \sl \fi\char36}}}{%
        \def\pdollar{\text{\ifdim \fontdimen\@ne\font >\z@
              \sl \else \normalshape \fi\char36}}}

\expandafter\ifx\csname fiverm\endcsname\relax
  \let\fiverm\fivrm
\fi
  
\let\!latexendpicture=\endpicture 
\let\!latexframe=\frame
\let\!latexlinethickness=\linethickness
\let\!latexmultiput=\multiput
\let\!latexput=\put
 
\def\@picture(#1,#2)(#3,#4){%
  \@picht #2\unitlength
  \setbox\@picbox\hbox to #1\unitlength\bgroup 
  \let\endpicture=\!latexendpicture
  \let\frame=\!latexframe
  \let\linethickness=\!latexlinethickness
  \let\multiput=\!latexmultiput
  \let\put=\!latexput
  \hskip -#3\unitlength \lower #4\unitlength \hbox\bgroup}

\catcode`@=12 \catcode`!=12

\catcode`!=11 
 
\def\PiC{P\kern-.12em\lower.5ex\hbox{I}\kern-.075emC}
\def\PiCTeX{\PiC\kern-.11em\TeX}

\def\!ifnextchar#1#2#3{%
  \let\!testchar=#1%
  \def\!first{#2}%
  \def\!second{#3}%
  \futurelet\!nextchar\!testnext}
\def\!testnext{%
  \ifx \!nextchar \!spacetoken 
    \let\!next=\!skipspacetestagain
  \else
    \ifx \!nextchar \!testchar
      \let\!next=\!first
    \else 
      \let\!next=\!second 
    \fi 
  \fi
  \!next}
\def\\{\!skipspacetestagain} 
  \expandafter\def\\ {\futurelet\!nextchar\!testnext} 
\def\\{\let\!spacetoken= } \\  

\def\!tfor#1:=#2\do#3{%
  \edef\!fortemp{#2}%
  \ifx\!fortemp\!empty 
    \else
    \!tforloop#2\!nil\!nil\!!#1{#3}%
  \fi}
\def\!tforloop#1#2\!!#3#4{%
  \def#3{#1}%
  \ifx #3\!nnil
    \let\!nextwhile=\!fornoop
  \else
    #4\relax
    \let\!nextwhile=\!tforloop
  \fi 
  \!nextwhile#2\!!#3{#4}}

\def\!etfor#1:=#2\do#3{%
  \def\!!tfor{\!tfor#1:=}%
  \edef\!!!tfor{#2}%
  \expandafter\!!tfor\!!!tfor\do{#3}}

\def\!cfor#1:=#2\do#3{%
  \edef\!fortemp{#2}%
  \ifx\!fortemp\!empty 
  \else
    \!cforloop#2,\!nil,\!nil\!!#1{#3}%
  \fi}
\def\!cforloop#1,#2\!!#3#4{%
  \def#3{#1}%
  \ifx #3\!nnil
    \let\!nextwhile=\!fornoop 
  \else
    #4\relax
    \let\!nextwhile=\!cforloop
  \fi
  \!nextwhile#2\!!#3{#4}}

\def\!ecfor#1:=#2\do#3{%
  \def\!!cfor{\!cfor#1:=}%
  \edef\!!!cfor{#2}%
  \expandafter\!!cfor\!!!cfor\do{#3}}

\def\!empty{}
\def\!nnil{\!nil}
\def\!fornoop#1\!!#2#3{}

\def\!ifempty#1#2#3{%
  \edef\!emptyarg{#1}%
  \ifx\!emptyarg\!empty
    #2%
  \else
    #3%
  \fi}
 
\def\!getnext#1\from#2{%
  \expandafter\!gnext#2\!#1#2}%
\def\!gnext\\#1#2\!#3#4{%
  \def#3{#1}%
  \def#4{#2\\{#1}}%
  \ignorespaces}

\def\!getnextvalueof#1\from#2{%
  \expandafter\!gnextv#2\!#1#2}%
\def\!gnextv\\#1#2\!#3#4{%
  #3=#1%
  \def#4{#2\\{#1}}%
  \ignorespaces}

\def\!copylist#1\to#2{%
  \expandafter\!!copylist#1\!#2}
\def\!!copylist#1\!#2{%
  \def#2{#1}\ignorespaces}

\def\!wlet#1=#2{%
  \let#1=#2 
  \wlog{\string#1=\string#2}}
 
\def\!listaddon#1#2{%
  \expandafter\!!listaddon#2\!{#1}#2}
\def\!!listaddon#1\!#2#3{%
  \def#3{#1\\#2}}
 
\def\!rightappend#1\withCS#2\to#3{\expandafter\!!rightappend#3\!#2{#1}#3}
\def\!!rightappend#1\!#2#3#4{\def#4{#1#2{#3}}}

\def\!leftappend#1\withCS#2\to#3{\expandafter\!!leftappend#3\!#2{#1}#3}
\def\!!leftappend#1\!#2#3#4{\def#4{#2{#3}#1}}

\def\!lop#1\to#2{\expandafter\!!lop#1\!#1#2}
\def\!!lop\\#1#2\!#3#4{\def#4{#1}\def#3{#2}}

\def\!loop#1\repeat{\def\!body{#1}\!iterate}
\def\!iterate{\!body\let\!next=\!iterate\else\let\!next=\relax\fi\!next}
 
\def\!!loop#1\repeat{\def\!!body{#1}\!!iterate}
\def\!!iterate{\!!body\let\!!next=\!!iterate\else\let\!!next=\relax\fi\!!next}
 
\def\!removept#1#2{\edef#2{\expandafter\!!removePT\the#1}}
{\catcode`p=12 \catcode`t=12 \gdef\!!removePT#1pt{#1}}

\def\placevalueinpts of <#1> in #2 {%
  \!removept{#1}{#2}}
 
\def\!mlap#1{\hbox to 0pt{\hss#1\hss}}
\def\!vmlap#1{\vbox to 0pt{\vss#1\vss}}
 
\def\!not#1{%
  #1\relax
    \!switchfalse
  \else
    \!switchtrue
  \fi
  \if!switch
  \ignorespaces}

\let\!!!wlog=\wlog              
\def\wlog#1{}    

\newdimen\headingtoplotskip     
\newdimen\linethickness         
\newdimen\longticklength        
\newdimen\plotsymbolspacing     
\newdimen\shortticklength       
\newdimen\stackleading          
\newdimen\tickstovaluesleading  
\newdimen\totalarclength        
\newdimen\valuestolabelleading  

\newbox\!boxA                   
\newbox\!boxB                   
\newbox\!picbox                 
\newbox\!plotsymbol             
\newbox\!putobject              
\newbox\!shadesymbol            

\newcount\!countA               
\newcount\!countB               
\newcount\!countC               
\newcount\!countD               
\newcount\!countE               
\newcount\!countF               
\newcount\!countG               
\newcount\!fiftypt              
\newcount\!intervalno           
\newcount\!npoints              
\newcount\!nsegments            
\newcount\!ntemp                
\newcount\!parity               
\newcount\!scalefactor          
\newcount\!tfs                  
\newcount\!tickcase             

\newdimen\!Xleft                
\newdimen\!Xright               
\newdimen\!Xsave                
\newdimen\!Ybot                 
\newdimen\!Ysave                
\newdimen\!Ytop                 
\newdimen\!angle                
\newdimen\!arclength            
\newdimen\!areabloc             
\newdimen\!arealloc             
\newdimen\!arearloc             
\newdimen\!areatloc             
\newdimen\!bshrinkage           
\newdimen\!checkbot             
\newdimen\!checkleft            
\newdimen\!checkright           
\newdimen\!checktop             
\newdimen\!dimenA               
\newdimen\!dimenB               
\newdimen\!dimenC               
\newdimen\!dimenD               
\newdimen\!dimenE               
\newdimen\!dimenF               
\newdimen\!dimenG               
\newdimen\!dimenH               
\newdimen\!dimenI               
\newdimen\!distacross           
\newdimen\!downlength           
\newdimen\!dp                   
\newdimen\!dshade               
\newdimen\!dxpos                
\newdimen\!dxprime              
\newdimen\!dypos                
\newdimen\!dyprime              
\newdimen\!ht                   
\newdimen\!leaderlength         
\newdimen\!lshrinkage           
\newdimen\!midarclength         
\newdimen\!offset               
\newdimen\!plotheadingoffset    
\newdimen\!plotsymbolxshift     
\newdimen\!plotsymbolyshift     
\newdimen\!plotxorigin          
\newdimen\!plotyorigin          
\newdimen\!rootten              
\newdimen\!rshrinkage           
\newdimen\!shadesymbolxshift    
\newdimen\!shadesymbolyshift    
\newdimen\!tenAa                
\newdimen\!tenAc                
\newdimen\!tenAe                
\newdimen\!tshrinkage           
\newdimen\!uplength             
\newdimen\!wd                   
\newdimen\!wmax                 
\newdimen\!wmin                 
\newdimen\!xB                   
\newdimen\!xC                   
\newdimen\!xE                   
\newdimen\!xM                   
\newdimen\!xS                   
\newdimen\!xaxislength          
\newdimen\!xdiff                
\newdimen\!xleft                
\newdimen\!xloc                 
\newdimen\!xorigin              
\newdimen\!xpivot               

\newdimen\!xpos                 
\newdimen\!xprime               
\newdimen\!xright               
\newdimen\!xshade               
\newdimen\!xshift               
\newdimen\!xtemp                
\newdimen\!xunit                
\newdimen\!xxE                  
\newdimen\!xxM                  
\newdimen\!xxS                  
\newdimen\!xxloc                
\newdimen\!yB                   
\newdimen\!yC                   
\newdimen\!yE                   
\newdimen\!yM                   
\newdimen\!yS                   
\newdimen\!yaxislength          
\newdimen\!ybot                 
\newdimen\!ydiff                
\newdimen\!yloc                 
\newdimen\!yorigin              
\newdimen\!ypivot               
\newdimen\!ypos                 
\newdimen\!yprime               
\newdimen\!yshade               
\newdimen\!yshift               
\newdimen\!ytemp                
\newdimen\!ytop                 
\newdimen\!yunit                
\newdimen\!yyE                  
\newdimen\!yyM                  
\newdimen\!yyS                  
\newdimen\!yyloc                
\newdimen\!zpt                  

\newif\if!axisvisible           
\newif\if!gridlinestoo          
\newif\if!keepPO                
\newif\if!placeaxislabel        
\newif\if!switch                
\newif\if!xswitch               

\newtoks\!axisLaBeL             
\newtoks\!keywordtoks           

\newwrite\!replotfile           

\newhelp\!keywordhelp{The keyword mentioned in the error message in unknown. 
Replace NEW KEYWORD in the indicated response by the keyword that 
should have been specified.}    

\!wlet\!!origin=\!xM                   
\!wlet\!!unit=\!uplength               
\!wlet\!Lresiduallength=\!dimenG       
\!wlet\!Rresiduallength=\!dimenF       
\!wlet\!axisLength=\!distacross        
\!wlet\!axisend=\!ydiff                
\!wlet\!axisstart=\!xdiff              
\!wlet\!axisxlevel=\!arclength         
\!wlet\!axisylevel=\!downlength        
\!wlet\!beta=\!dimenE                  
\!wlet\!gamma=\!dimenF                 
\!wlet\!shadexorigin=\!plotxorigin     
\!wlet\!shadeyorigin=\!plotyorigin     
\!wlet\!ticklength=\!xS                
\!wlet\!ticklocation=\!xE              
\!wlet\!ticklocationincr=\!yE          
\!wlet\!tickwidth=\!yS                 
\!wlet\!totalleaderlength=\!dimenE     
\!wlet\!xone=\!xprime                  
\!wlet\!xtwo=\!dxprime                 
\!wlet\!ySsave=\!yM                    
\!wlet\!ybB=\!yB                       
\!wlet\!ybC=\!yC                       
\!wlet\!ybE=\!yE                       
\!wlet\!ybM=\!yM                       
\!wlet\!ybS=\!yS                       
\!wlet\!ybpos=\!yyloc                  
\!wlet\!yone=\!yprime                  
\!wlet\!ytB=\!xB                       
\!wlet\!ytC=\!xC                       
\!wlet\!ytE=\!downlength               
\!wlet\!ytM=\!arclength                
\!wlet\!ytS=\!distacross               
\!wlet\!ytpos=\!xxloc                  
\!wlet\!ytwo=\!dyprime                 

\!zpt=0pt                              
\!xunit=1pt
\!yunit=1pt
\!arearloc=\!xunit
\!areatloc=\!yunit
\!dshade=5pt
\!leaderlength=24in
\!tfs=256                              
\!wmax=5.3pt                           
\!wmin=2.7pt                           
\!xaxislength=\!xunit
\!xpivot=\!zpt
\!yaxislength=\!yunit 
\!ypivot=\!zpt
  \!dimenA=50pt \!fiftypt=\!dimenA     

\!rootten=3.162278pt                   
\!tenAa=8.690286pt                     
\!tenAc=2.773839pt                     
\!tenAe=2.543275pt                     

\def\!cosrotationangle{1}      
\def\!sinrotationangle{0}      
\def\!xpivotcoord{0}           
\def\!xref{0}                  
\def\!xshadesave{0}            
\def\!ypivotcoord{0}           
\def\!yref{0}                  
\def\!yshadesave{0}            
\def\!zero{0}                  

\let\wlog=\!!!wlog
\def\normalgraphs{%
  \longticklength=.4\baselineskip
  \shortticklength=.25\baselineskip
  \tickstovaluesleading=.25\baselineskip
  \valuestolabelleading=.8\baselineskip
  \linethickness=.4pt
  \stackleading=.17\baselineskip
  \headingtoplotskip=1.5\baselineskip
  \visibleaxes
  \ticksout
  \nogridlines
  \unloggedticks}
\def\setplotarea x from #1 to #2, y from #3 to #4 {%
  \!arealloc=\!M{#1}\!xunit \advance \!arealloc -\!xorigin
  \!areabloc=\!M{#3}\!yunit \advance \!areabloc -\!yorigin
  \!arearloc=\!M{#2}\!xunit \advance \!arearloc -\!xorigin
  \!areatloc=\!M{#4}\!yunit \advance \!areatloc -\!yorigin
  \!initinboundscheck
  \!xaxislength=\!arearloc  \advance\!xaxislength -\!arealloc
  \!yaxislength=\!areatloc  \advance\!yaxislength -\!areabloc
  \!plotheadingoffset=\!zpt
  \!dimenput {{\setbox0=\hbox{}\wd0=\!xaxislength\ht0=\!yaxislength\box0}}
     [bl] (\!arealloc,\!areabloc)}
\def\visibleaxes{%
  \def\!axisvisibility{\!axisvisibletrue}}

\def\!fixkeyword#1{%
  \errhelp=\!keywordhelp
  \errmessage{Unrecognized keyword `#1': \the\!keywordtoks{NEW KEYWORD}'}}

\!keywordtoks={enter `i\fixkeyword}

\def\fixkeyword#1{%
  \!nextkeyword#1 }

\def\axis {%
  \def\!nextkeyword##1 {%
    \expandafter\ifx\csname !axis##1\endcsname \relax
      \def\!next{\!fixkeyword{##1}}%
    \else
      \def\!next{\csname !axis##1\endcsname}%
    \fi
    \!next}%
  \!offset=\!zpt
  \!axisvisibility
  \!placeaxislabelfalse
  \!nextkeyword}

\def\!axisbottom{%
  \!axisylevel=\!areabloc
  \def\!tickxsign{0}%
  \def\!tickysign{-}%
  \def\!axissetup{\!axisxsetup}%
  \def\!axislabeltbrl{t}%
  \!nextkeyword}

\def\!axistop{%
  \!axisylevel=\!areatloc
  \def\!tickxsign{0}%
  \def\!tickysign{+}%
  \def\!axissetup{\!axisxsetup}%
  \def\!axislabeltbrl{b}%
  \!nextkeyword}

\def\!axisleft{%
  \!axisxlevel=\!arealloc
  \def\!tickxsign{-}%
  \def\!tickysign{0}%
  \def\!axissetup{\!axisysetup}%
  \def\!axislabeltbrl{r}%
  \!nextkeyword}

\def\!axisright{%
  \!axisxlevel=\!arearloc
  \def\!tickxsign{+}%
  \def\!tickysign{0}%
  \def\!axissetup{\!axisysetup}%
  \def\!axislabeltbrl{l}%
  \!nextkeyword}

\def\!axisshiftedto#1=#2 {%
  \if 0\!tickxsign
    \!axisylevel=\!M{#2}\!yunit
    \advance\!axisylevel -\!yorigin
  \else
    \!axisxlevel=\!M{#2}\!xunit
    \advance\!axisxlevel -\!xorigin
  \fi
  \!nextkeyword}

\def\!axisvisible{%
  \!axisvisibletrue  
  \!nextkeyword}

\def\!axisinvisible{%
  \!axisvisiblefalse
  \!nextkeyword}

\def\!axislabel#1 {%
  \!axisLaBeL={#1}%
  \!placeaxislabeltrue
  \!nextkeyword}

\expandafter\def\csname !axis/\endcsname{%
  \!axissetup 
  \if!placeaxislabel
    \!placeaxislabel
  \fi
  \if +\!tickysign 
    \!dimenA=\!axisylevel
    \advance\!dimenA \!offset 
    \advance\!dimenA -\!areatloc 
    \ifdim \!dimenA>\!plotheadingoffset
      \!plotheadingoffset=\!dimenA 
    \fi
  \fi}

\def\grid #1 #2 {%
  \!countA=#1\advance\!countA 1
  \axis bottom invisible ticks length <\!zpt> andacross quantity {\!countA} /
  \!countA=#2\advance\!countA 1
  \axis left   invisible ticks length <\!zpt> andacross quantity {\!countA} / }

\def\plotheading#1 {%
  \advance\!plotheadingoffset \headingtoplotskip
  \!dimenput {#1} [B] <.5\!xaxislength,\!plotheadingoffset>
    (\!arealloc,\!areatloc)}

\def\!axisxsetup{%
  \!axisxlevel=\!arealloc
  \!axisstart=\!arealloc
  \!axisend=\!arearloc
  \!axisLength=\!xaxislength
  \!!origin=\!xorigin
  \!!unit=\!xunit
  \!xswitchtrue
  \if!axisvisible 
    \!makeaxis
  \fi}

\def\!axisysetup{%
  \!axisylevel=\!areabloc
  \!axisstart=\!areabloc
  \!axisend=\!areatloc
  \!axisLength=\!yaxislength
  \!!origin=\!yorigin
  \!!unit=\!yunit
  \!xswitchfalse
  \if!axisvisible
    \!makeaxis
  \fi}

\def\!makeaxis{%
  \setbox\!boxA=\hbox{
    \beginpicture
      \!setdimenmode
      \setcoordinatesystem point at {\!zpt} {\!zpt}   
      \putrule from {\!zpt} {\!zpt} to
        {\!tickysign\!tickysign\!axisLength} 
        {\!tickxsign\!tickxsign\!axisLength}
    \endpicturesave <\!Xsave,\!Ysave>}%
    \wd\!boxA=\!zpt
    \!placetick\!axisstart}

\def\!placeaxislabel{%
  \advance\!offset \valuestolabelleading
  \if!xswitch
    \!dimenput {\the\!axisLaBeL} [\!axislabeltbrl]
      <.5\!axisLength,\!tickysign\!offset> (\!axisxlevel,\!axisylevel)
    \advance\!offset \!dp  
    \advance\!offset \!ht  
  \else
    \!dimenput {\the\!axisLaBeL} [\!axislabeltbrl]
      <\!tickxsign\!offset,.5\!axisLength> (\!axisxlevel,\!axisylevel)
  \fi
  \!axisLaBeL={}}

\def\arrow <#1> [#2,#3]{%
  \!ifnextchar<{\!arrow{#1}{#2}{#3}}{\!arrow{#1}{#2}{#3}<\!zpt,\!zpt> }}

\def\!arrow#1#2#3<#4,#5> from #6 #7 to #8 #9 {%
  \!xloc=\!M{#8}\!xunit   
  \!yloc=\!M{#9}\!yunit
  \!dxpos=\!xloc  \!dimenA=\!M{#6}\!xunit  \advance \!dxpos -\!dimenA
  \!dypos=\!yloc  \!dimenA=\!M{#7}\!yunit  \advance \!dypos -\!dimenA
  \let\!MAH=\!M
  \!setdimenmode
  \!xshift=#4\relax  \!yshift=#5\relax
  \!reverserotateonly\!xshift\!yshift
  \advance\!xshift\!xloc  \advance\!yshift\!yloc
  \!xS=-\!dxpos  \advance\!xS\!xshift
  \!yS=-\!dypos  \advance\!yS\!yshift
  \!start (\!xS,\!yS)
  \!ljoin (\!xshift,\!yshift)
  \!Pythag\!dxpos\!dypos\!arclength
  \!divide\!dxpos\!arclength\!dxpos  
  \!dxpos=32\!dxpos  \!removept\!dxpos\!!cos
  \!divide\!dypos\!arclength\!dypos  
  \!dypos=32\!dypos  \!removept\!dypos\!!sin
  \!halfhead{#1}{#2}{#3}
  \!halfhead{#1}{-#2}{-#3}
  \let\!M=\!MAH
  \ignorespaces}
  \def\!halfhead#1#2#3{%
    \!dimenC=-#1%
    \divide \!dimenC 2 
    \!dimenD=#2\!dimenC
    \!rotate(\!dimenC,\!dimenD)by(\!!cos,\!!sin)to(\!xM,\!yM)
    \!dimenC=-#1
    \!dimenD=#3\!dimenC
    \!dimenD=.5\!dimenD
    \!rotate(\!dimenC,\!dimenD)by(\!!cos,\!!sin)to(\!xE,\!yE)
    \!start (\!xshift,\!yshift)
    \advance\!xM\!xshift  \advance\!yM\!yshift
    \advance\!xE\!xshift  \advance\!yE\!yshift
    \!qjoin (\!xM,\!yM) (\!xE,\!yE) 
    \ignorespaces}

\def\betweenarrows #1#2 from #3 #4 to #5 #6 {%
  \!xloc=\!M{#3}\!xunit  \!xxloc=\!M{#5}\!xunit%
  \!yloc=\!M{#4}\!yunit  \!yyloc=\!M{#6}\!yunit%
  \!dxpos=\!xxloc  \advance\!dxpos by -\!xloc
  \!dypos=\!yyloc  \advance\!dypos by -\!yloc
  \advance\!xloc .5\!dxpos
  \advance\!yloc .5\!dypos
  \let\!MBA=\!M
  \!setdimenmode
  \ifdim\!dypos=\!zpt
    \ifdim\!dxpos<\!zpt \!dxpos=-\!dxpos \fi
    \put {\!lrarrows{\!dxpos}{#1}}#2{} at {\!xloc} {\!yloc}
  \else
    \ifdim\!dxpos=\!zpt
      \ifdim\!dypos<\!zpt \!dypos=-\!zpt \fi
      \put {\!udarrows{\!dypos}{#1}}#2{} at {\!xloc} {\!yloc}
    \fi
  \fi
  \let\!M=\!MBA
  \ignorespaces}

\def\!lrarrows#1#2{
  {\setbox\!boxA=\hbox{$\mkern-2mu\mathord-\mkern-2mu$}%
   \setbox\!boxB=\hbox{$\leftarrow$}\!dimenE=\ht\!boxB
   \setbox\!boxB=\hbox{}\ht\!boxB=2\!dimenE
   \hbox to #1{$\mathord\leftarrow\mkern-6mu
     \cleaders\copy\!boxA\hfil
     \mkern-6mu\mathord-$%
     \kern.4em $\vcenter{\box\!boxB}$$\vcenter{\hbox{#2}}$\kern.4em
     $\mathord-\mkern-6mu
     \cleaders\copy\!boxA\hfil
     \mkern-6mu\mathord\rightarrow$}}}

\def\!udarrows#1#2{
  {\setbox\!boxB=\hbox{#2}%
   \setbox\!boxA=\hbox to \wd\!boxB{\hss$\vert$\hss}%
   \!dimenE=\ht\!boxA \advance\!dimenE \dp\!boxA \divide\!dimenE 2
   \vbox to #1{\offinterlineskip
      \vskip .05556\!dimenE
      \hbox to \wd\!boxB{\hss$\mkern.4mu\uparrow$\hss}\vskip-\!dimenE
      \cleaders\copy\!boxA\vfil
      \vskip-\!dimenE\copy\!boxA
      \vskip\!dimenE\copy\!boxB\vskip.4em
      \copy\!boxA\vskip-\!dimenE
      \cleaders\copy\!boxA\vfil
      \vskip-\!dimenE \hbox to \wd\!boxB{\hss$\mkern.4mu\downarrow$\hss}
      \vskip .05556\!dimenE}}}

\def\putbar#1breadth <#2> from #3 #4 to #5 #6 {%
  \!xloc=\!M{#3}\!xunit  \!xxloc=\!M{#5}\!xunit%
  \!yloc=\!M{#4}\!yunit  \!yyloc=\!M{#6}\!yunit%
  \!dypos=\!yyloc  \advance\!dypos by -\!yloc
  \!dimenI=#2  
  \ifdim \!dimenI=\!zpt 
    \putrule#1from {#3} {#4} to {#5} {#6} 
  \else 
    \let\!MBar=\!M
    \!setdimenmode 
    \divide\!dimenI 2
    \ifdim \!dypos=\!zpt             
      \advance \!yloc -\!dimenI 
      \advance \!yyloc \!dimenI
    \else
      \advance \!xloc -\!dimenI 
      \advance \!xxloc \!dimenI
    \fi
    \putrectangle#1corners at {\!xloc} {\!yloc} and {\!xxloc} {\!yyloc}
    \let\!M=\!MBar 
  \fi
  \ignorespaces}

\def\setbars#1breadth <#2> baseline at #3 = #4 {%
  \edef\!barshift{#1}%
  \edef\!barbreadth{#2}%
  \edef\!barorientation{#3}%
  \edef\!barbaseline{#4}%
  \def\!bardobaselabel{\!bardoendlabel}%
  \def\!bardoendlabel{\!barfinish}%
  \let\!drawcurve=\!barcurve
  \!setbars}
\def\!setbars{%
  \futurelet\!nextchar\!!setbars}
\def\!!setbars{%
  \if b\!nextchar
    \def\!!!setbars{\!setbarsbget}%
  \else 
    \if e\!nextchar
      \def\!!!setbars{\!setbarseget}%
    \else
      \def\!!!setbars{\relax}%
    \fi
  \fi
  \!!!setbars}
\def\!setbarsbget baselabels (#1) {%
  \def\!barbaselabelorientation{#1}%
  \def\!bardobaselabel{\!!bardobaselabel}%
  \!setbars}
\def\!setbarseget endlabels (#1) {%
  \edef\!barendlabelorientation{#1}%
  \def\!bardoendlabel{\!!bardoendlabel}%
  \!setbars}

\def\!barcurve #1 #2 {%
  \if y\!barorientation
    \def\!basexarg{#1}%
    \def\!baseyarg{\!barbaseline}%
  \else
    \def\!basexarg{\!barbaseline}%
    \def\!baseyarg{#2}%
  \fi
  \expandafter\putbar\!barshift breadth <\!barbreadth> from {\!basexarg}
    {\!baseyarg} to {#1} {#2}
  \def\!endxarg{#1}%
  \def\!endyarg{#2}%
  \!bardobaselabel}

\def\!!bardobaselabel "#1" {%
  \put {#1}\!barbaselabelorientation{} at {\!basexarg} {\!baseyarg}
  \!bardoendlabel}
 
\def\!!bardoendlabel "#1" {%
  \put {#1}\!barendlabelorientation{} at {\!endxarg} {\!endyarg}
  \!barfinish}

\def\!barfinish{%
  \!ifnextchar/{\!finish}{\!barcurve}}

\def\putrectangle{%
  \!ifnextchar<{\!putrectangle}{\!putrectangle<\!zpt,\!zpt> }}
\def\!putrectangle<#1,#2> corners at #3 #4 and #5 #6 {%
  \!xone=\!M{#3}\!xunit  \!xtwo=\!M{#5}\!xunit%
  \!yone=\!M{#4}\!yunit  \!ytwo=\!M{#6}\!yunit%
  \ifdim \!xtwo<\!xone
    \!dimenI=\!xone  \!xone=\!xtwo  \!xtwo=\!dimenI
  \fi
  \ifdim \!ytwo<\!yone
    \!dimenI=\!yone  \!yone=\!ytwo  \!ytwo=\!dimenI
  \fi
  \!dimenI=#1\relax  \advance\!xone\!dimenI  \advance\!xtwo\!dimenI
  \!dimenI=#2\relax  \advance\!yone\!dimenI  \advance\!ytwo\!dimenI
  \let\!MRect=\!M
  \!setdimenmode
  \!shaderectangle
  \!dimenI=.5\linethickness
  \advance \!xone  -\!dimenI
  \advance \!xtwo   \!dimenI
  \putrule from {\!xone} {\!yone} to {\!xtwo} {\!yone} 
  \putrule from {\!xone} {\!ytwo} to {\!xtwo} {\!ytwo} 
  \advance \!xone   \!dimenI
  \advance \!xtwo  -\!dimenI%
  \advance \!yone  -\!dimenI
  \advance \!ytwo   \!dimenI
  \putrule from {\!xone} {\!yone} to {\!xone} {\!ytwo} 
  \putrule from {\!xtwo} {\!yone} to {\!xtwo} {\!ytwo} 
  \let\!M=\!MRect
  \ignorespaces}
 
\def\shaderectanglesoff{%
  \def\!shaderectangle{}%
  \ignorespaces}

\shaderectanglesoff
 
\def\!!shaderectangle{%
  \!dimenA=\!xtwo  \advance \!dimenA -\!xone
  \!dimenB=\!ytwo  \advance \!dimenB -\!yone
  \ifdim \!dimenA<\!dimenB
    \!startvshade (\!xone,\!yone,\!ytwo)
    \!lshade      (\!xtwo,\!yone,\!ytwo)
  \else
    \!starthshade (\!yone,\!xone,\!xtwo)
    \!lshade      (\!ytwo,\!xone,\!xtwo)
  \fi
  \ignorespaces}
  
\def\frame{%
  \!ifnextchar<{\!frame}{\!frame<\!zpt> }}
\long\def\!frame<#1> #2{%
  \beginpicture
    \setcoordinatesystem units <1pt,1pt> point at 0 0 
    \put {#2} [Bl] at 0 0 
    \!dimenA=#1\relax
    \!dimenB=\!wd \advance \!dimenB \!dimenA
    \!dimenC=\!ht \advance \!dimenC \!dimenA
    \!dimenD=\!dp \advance \!dimenD \!dimenA
    \let\!MFr=\!M
    \!setdimenmode
    \putrectangle corners at {-\!dimenA} {-\!dimenD} and {\!dimenB} {\!dimenC}
    \!setcoordmode
    \let\!M=\!MFr
  \endpicture
  \ignorespaces}
 
\def\rectangle <#1> <#2> {%
  \setbox0=\hbox{}\wd0=#1\ht0=#2\frame {\box0}}

\def\!plotfromfile"#1"{%
  \expandafter\!drawcurve \input #1 /}

\def\setquadratic{%
  \let\!drawcurve=\!qcurve
  \let\!!Shade=\!!qShade
  \let\!!!Shade=\!!!qShade}

\def\setlinear{%
  \let\!drawcurve=\!lcurve
  \let\!!Shade=\!!lShade
  \let\!!!Shade=\!!!lShade}

\def\sethistograms{%
  \let\!drawcurve=\!hcurve}

\def\!qcurve #1 #2 {%
  \!start (#1,#2)
  \!Qjoin}
\def\!Qjoin#1 #2 #3 #4 {%
  \!qjoin (#1,#2) (#3,#4)             
  \!ifnextchar/{\!finish}{\!Qjoin}}

\def\!lcurve #1 #2 {%
  \!start (#1,#2)
  \!Ljoin}
\def\!Ljoin#1 #2 {%
  \!ljoin (#1,#2)                    
  \!ifnextchar/{\!finish}{\!Ljoin}}

\def\!finish/{\ignorespaces}

\def\!hcurve #1 #2 {%
  \edef\!hxS{#1}%
  \edef\!hyS{#2}%
  \!hjoin}
\def\!hjoin#1 #2 {%
  \putrectangle corners at {\!hxS} {\!hyS} and {#1} {#2}
  \edef\!hxS{#1}%
  \!ifnextchar/{\!finish}{\!hjoin}}

\def\vshade #1 #2 #3 {%
  \!startvshade (#1,#2,#3)
  \!Shadewhat}

\def\hshade #1 #2 #3 {%
  \!starthshade (#1,#2,#3)
  \!Shadewhat}

\def\!Shadewhat{%
  \futurelet\!nextchar\!Shade}
\def\!Shade{%
  \if <\!nextchar
    \def\!nextShade{\!!Shade}%
  \else
    \if /\!nextchar
      \def\!nextShade{\!finish}%
    \else
      \def\!nextShade{\!!!Shade}%
    \fi
  \fi
  \!nextShade}
\def\!!lShade<#1> #2 #3 #4 {%
  \!lshade <#1> (#2,#3,#4)                 
  \!Shadewhat}
\def\!!!lShade#1 #2 #3 {%
  \!lshade (#1,#2,#3)
  \!Shadewhat} 
\def\!!qShade<#1> #2 #3 #4 #5 #6 #7 {%
  \!qshade <#1> (#2,#3,#4) (#5,#6,#7)      
  \!Shadewhat}
\def\!!!qShade#1 #2 #3 #4 #5 #6 {%
  \!qshade (#1,#2,#3) (#4,#5,#6)
  \!Shadewhat} 

\setlinear

\def\setdashpattern <#1>{%
  \def\!Flist{}\def\!Blist{}\def\!UDlist{}%
  \!countA=0
  \!ecfor\!item:=#1\do{%
    \!dimenA=\!item\relax
    \expandafter\!rightappend\the\!dimenA\withCS{\\}\to\!UDlist%
    \advance\!countA  1
    \ifodd\!countA
      \expandafter\!rightappend\the\!dimenA\withCS{\!Rule}\to\!Flist%
      \expandafter\!leftappend\the\!dimenA\withCS{\!Rule}\to\!Blist%
    \else 
      \expandafter\!rightappend\the\!dimenA\withCS{\!Skip}\to\!Flist%
      \expandafter\!leftappend\the\!dimenA\withCS{\!Skip}\to\!Blist%
    \fi}%
  \!leaderlength=\!zpt
  \def\!Rule##1{\advance\!leaderlength  ##1}%
  \def\!Skip##1{\advance\!leaderlength  ##1}%
  \!Flist%
  \ifdim\!leaderlength>\!zpt 
  \else
    \def\!Flist{\!Skip{24in}}\def\!Blist{\!Skip{24in}}\ignorespaces
    \def\!UDlist{\\{\!zpt}\\{24in}}\ignorespaces
    \!leaderlength=24in
  \fi
  \!dashingon}

\def\!dashingon{%
  \def\!advancedashing{\!!advancedashing}%
  \def\!drawlinearsegment{\!lineardashed}%
  \def\!puthline{\!putdashedhline}%
  \def\!putvline{\!putdashedvline}%
  \ignorespaces}%
\def\!dashingoff{%
  \def\!advancedashing{\relax}%
  \def\!drawlinearsegment{\!linearsolid}%
  \def\!puthline{\!putsolidhline}%
  \def\!putvline{\!putsolidvline}%
  \ignorespaces}

\def\setdots{%
  \!ifnextchar<{\!setdots}{\!setdots<5pt>}}
\def\!setdots<#1>{%
  \!dimenB=#1\advance\!dimenB -\plotsymbolspacing
  \ifdim\!dimenB<\!zpt
    \!dimenB=\!zpt
  \fi
\setdashpattern <\plotsymbolspacing,\!dimenB>}
 
\def\setdotsnear <#1> for <#2>{%
  \!dimenB=#2\relax  \advance\!dimenB -.05pt  
  \!dimenC=#1\relax  \!countA=\!dimenC 
  \!dimenD=\!dimenB  \advance\!dimenD .5\!dimenC  \!countB=\!dimenD
  \divide \!countB  \!countA
  \ifnum 1>\!countB 
    \!countB=1
  \fi
  \divide\!dimenB  \!countB
  \setdots <\!dimenB>}
 
\def\setdashes{%
  \!ifnextchar<{\!setdashes}{\!setdashes<5pt>}}
\def\!setdashes<#1>{\setdashpattern <#1,#1>}
 
\def\setdashesnear <#1> for <#2>{%
  \!dimenB=#2\relax  
  \!dimenC=#1\relax  \!countA=\!dimenC 
  \!dimenD=\!dimenB  \advance\!dimenD .5\!dimenC  \!countB=\!dimenD
  \divide \!countB  \!countA
  \ifodd \!countB 
  \else 
    \advance \!countB  1
  \fi
  \divide\!dimenB  \!countB
  \setdashes <\!dimenB>}
 
\def\setsolid{%
  \def\!Flist{\!Rule{24in}}\def\!Blist{\!Rule{24in}}%
  \def\!UDlist{\\{24in}\\{\!zpt}}%
  \!dashingoff}  
\setsolid

\def\!divide#1#2#3{%
  \!dimenB=#1
  \!dimenC=#2
  \!dimenD=\!dimenB
  \divide \!dimenD \!dimenC
  \!dimenA=\!dimenD
  \multiply\!dimenD \!dimenC
  \advance\!dimenB -\!dimenD
  \!dimenD=\!dimenC
    \ifdim\!dimenD<\!zpt \!dimenD=-\!dimenD 
  \fi
  \ifdim\!dimenD<64pt
    \!divstep[\!tfs]\!divstep[\!tfs]%
  \else 
    \!!divide
  \fi
  #3=\!dimenA\ignorespaces}

\def\!!divide{%
  \ifdim\!dimenD<256pt
    \!divstep[64]\!divstep[32]\!divstep[32]%
  \else 
    \!divstep[8]\!divstep[8]\!divstep[8]\!divstep[8]\!divstep[8]%
    \!dimenA=2\!dimenA
  \fi}

\def\!divstep[#1]{
  \!dimenB=#1\!dimenB
  \!dimenD=\!dimenB
    \divide \!dimenD by \!dimenC
  \!dimenA=#1\!dimenA
    \advance\!dimenA by \!dimenD%
  \multiply\!dimenD by \!dimenC
    \advance\!dimenB by -\!dimenD}
 
\def\Divide <#1> by <#2> forming <#3> {%
  \!divide{#1}{#2}{#3}}

\def\ellipticalarc axes ratio #1:#2 #3 degrees from #4 #5 center at #6 #7 {%
  \!angle=#3pt\relax
  \ifdim\!angle>\!zpt 
    \def\!sign{}
  \else 
    \def\!sign{-}\!angle=-\!angle
  \fi
  \!xxloc=\!M{#6}\!xunit
  \!yyloc=\!M{#7}\!yunit     
  \!xxS=\!M{#4}\!xunit
  \!yyS=\!M{#5}\!yunit
  \advance\!xxS -\!xxloc
  \advance\!yyS -\!yyloc
  \!divide\!xxS{#1pt}\!xxS 
  \!divide\!yyS{#2pt}\!yyS 
  \let\!MC=\!M
  \!setdimenmode
  \!xS=#1\!xxS  \advance\!xS\!xxloc
  \!yS=#2\!yyS  \advance\!yS\!yyloc
  \!start (\!xS,\!yS)%
  \!loop\ifdim\!angle>14.9999pt
    \!rotate(\!xxS,\!yyS)by(\!cos,\!sign\!sin)to(\!xxM,\!yyM) 
    \!rotate(\!xxM,\!yyM)by(\!cos,\!sign\!sin)to(\!xxE,\!yyE)
    \!xM=#1\!xxM  \advance\!xM\!xxloc  \!yM=#2\!yyM  \advance\!yM\!yyloc
    \!xE=#1\!xxE  \advance\!xE\!xxloc  \!yE=#2\!yyE  \advance\!yE\!yyloc
    \!qjoin (\!xM,\!yM) (\!xE,\!yE)
    \!xxS=\!xxE  \!yyS=\!yyE 
    \advance \!angle -15pt
  \repeat
  \ifdim\!angle>\!zpt
    \!angle=100.53096\!angle
    \divide \!angle 360 
    \!sinandcos\!angle\!!sin\!!cos
    \!rotate(\!xxS,\!yyS)by(\!!cos,\!sign\!!sin)to(\!xxM,\!yyM) 
    \!rotate(\!xxM,\!yyM)by(\!!cos,\!sign\!!sin)to(\!xxE,\!yyE)
    \!xM=#1\!xxM  \advance\!xM\!xxloc  \!yM=#2\!yyM  \advance\!yM\!yyloc
    \!xE=#1\!xxE  \advance\!xE\!xxloc  \!yE=#2\!yyE  \advance\!yE\!yyloc
    \!qjoin (\!xM,\!yM) (\!xE,\!yE)
  \fi
  \let\!M=\!MC
  \ignorespaces}

\def\!rotate(#1,#2)by(#3,#4)to(#5,#6){%
  \!dimenA=#3#1\advance \!dimenA -#4#2
  \!dimenB=#3#2\advance \!dimenB  #4#1
  \divide \!dimenA 32  \divide \!dimenB 32 
  #5=\!dimenA  #6=\!dimenB
  \ignorespaces}
\def\!sin{4.17684}
\def\!cos{31.72624}

\def\!sinandcos#1#2#3{%
 \!dimenD=#1
 \!dimenA=\!dimenD
 \!dimenB=32pt
 \!removept\!dimenD\!value
 \!dimenC=\!dimenD
 \!dimenC=\!value\!dimenC \divide\!dimenC by 64 
 \advance\!dimenB by -\!dimenC
 \!dimenC=\!value\!dimenC \divide\!dimenC by 96 
 \advance\!dimenA by -\!dimenC
 \!dimenC=\!value\!dimenC \divide\!dimenC by 128 
 \advance\!dimenB by \!dimenC%
 \!removept\!dimenA#2
 \!removept\!dimenB#3
 \ignorespaces}

\def\putrule#1from #2 #3 to #4 #5 {%
  \!xloc=\!M{#2}\!xunit  \!xxloc=\!M{#4}\!xunit%
  \!yloc=\!M{#3}\!yunit  \!yyloc=\!M{#5}\!yunit%
  \!dxpos=\!xxloc  \advance\!dxpos by -\!xloc
  \!dypos=\!yyloc  \advance\!dypos by -\!yloc
  \ifdim\!dypos=\!zpt
    \def\!!Line{\!puthline{#1}}\ignorespaces
  \else
    \ifdim\!dxpos=\!zpt
      \def\!!Line{\!putvline{#1}}\ignorespaces
    \else 
       \def\!!Line{}
    \fi
  \fi
  \let\!ML=\!M
  \!setdimenmode
  \!!Line%
  \let\!M=\!ML
  \ignorespaces}

\def\!putsolidhline#1{%
  \ifdim\!dxpos>\!zpt 
    \put{\!hline\!dxpos}#1[l] at {\!xloc} {\!yloc}
  \else 
    \put{\!hline{-\!dxpos}}#1[l] at {\!xxloc} {\!yyloc}
  \fi
  \ignorespaces}
 
\def\!putsolidvline#1{%
  \ifdim\!dypos>\!zpt 
    \put{\!vline\!dypos}#1[b] at {\!xloc} {\!yloc}
  \else 
    \put{\!vline{-\!dypos}}#1[b] at {\!xxloc} {\!yyloc}
  \fi
  \ignorespaces}
 
\def\!hline#1{\hbox to #1{\leaders \hrule height\linethickness\hfill}}
\def\!vline#1{\vbox to #1{\leaders \vrule width\linethickness\vfill}}

\def\!putdashedhline#1{%
  \ifdim\!dxpos>\!zpt 
    \!DLsetup\!Flist\!dxpos
    \put{\hbox to \!totalleaderlength{\!hleaders}\!hpartialpattern\!Rtrunc}
      #1[l] at {\!xloc} {\!yloc} 
  \else 
    \!DLsetup\!Blist{-\!dxpos}
    \put{\!hpartialpattern\!Ltrunc\hbox to \!totalleaderlength{\!hleaders}}
      #1[r] at {\!xloc} {\!yloc} 
  \fi
  \ignorespaces}
 
\def\!putdashedvline#1{%
  \!dypos=-\!dypos
  \ifdim\!dypos>\!zpt 
    \!DLsetup\!Flist\!dypos 
    \put{\vbox{\vbox to \!totalleaderlength{\!vleaders}
      \!vpartialpattern\!Rtrunc}}#1[t] at {\!xloc} {\!yloc} 
  \else 
    \!DLsetup\!Blist{-\!dypos}
    \put{\vbox{\!vpartialpattern\!Ltrunc
      \vbox to \!totalleaderlength{\!vleaders}}}#1[b] at {\!xloc} {\!yloc} 
  \fi
  \ignorespaces}

\def\!DLsetup#1#2{
  \let\!RSlist=#1
  \!countB=#2
  \!countA=\!leaderlength
  \divide\!countB by \!countA
  \!totalleaderlength=\!countB\!leaderlength
  \!Rresiduallength=#2%
  \advance \!Rresiduallength by -\!totalleaderlength
  \!Lresiduallength=\!leaderlength
  \advance \!Lresiduallength by -\!Rresiduallength
  \ignorespaces}
 
\def\!hleaders{%
  \def\!Rule##1{\vrule height\linethickness width##1}%
  \def\!Skip##1{\hskip##1}%
  \leaders\hbox{\!RSlist}\hfill}
 
\def\!hpartialpattern#1{%
  \!dimenA=\!zpt \!dimenB=\!zpt 
  \def\!Rule##1{#1{##1}\vrule height\linethickness width\!dimenD}%
  \def\!Skip##1{#1{##1}\hskip\!dimenD}%
  \!RSlist}
 
\def\!vleaders{%
  \def\!Rule##1{\hrule width\linethickness height##1}%
  \def\!Skip##1{\vskip##1}%
  \leaders\vbox{\!RSlist}\vfill}
 
\def\!vpartialpattern#1{%
  \!dimenA=\!zpt \!dimenB=\!zpt 
  \def\!Rule##1{#1{##1}\hrule width\linethickness height\!dimenD}%
  \def\!Skip##1{#1{##1}\vskip\!dimenD}%
  \!RSlist}
 
\def\!Rtrunc#1{\!trunc{#1}>\!Rresiduallength}
\def\!Ltrunc#1{\!trunc{#1}<\!Lresiduallength}
 
\def\!trunc#1#2#3{%
  \!dimenA=\!dimenB         
  \advance\!dimenB by #1%
  \!dimenD=\!dimenB  \ifdim\!dimenD#2#3\!dimenD=#3\fi
  \!dimenC=\!dimenA  \ifdim\!dimenC#2#3\!dimenC=#3\fi
  \advance \!dimenD by -\!dimenC}

\def\!start (#1,#2){%
  \!plotxorigin=\!xorigin  \advance \!plotxorigin by \!plotsymbolxshift
  \!plotyorigin=\!yorigin  \advance \!plotyorigin by \!plotsymbolyshift
  \!xS=\!M{#1}\!xunit \!yS=\!M{#2}\!yunit
  \!rotateaboutpivot\!xS\!yS
  \!copylist\!UDlist\to\!!UDlist
  \!getnextvalueof\!downlength\from\!!UDlist
  \!distacross=\!zpt
  \!intervalno=0 
  \global\totalarclength=\!zpt
  \ignorespaces}

\def\!ljoin (#1,#2){%
  \advance\!intervalno by 1
  \!xE=\!M{#1}\!xunit \!yE=\!M{#2}\!yunit
  \!rotateaboutpivot\!xE\!yE
  \!xdiff=\!xE \advance \!xdiff by -\!xS
  \!ydiff=\!yE \advance \!ydiff by -\!yS
  \!Pythag\!xdiff\!ydiff\!arclength
  \global\advance \totalarclength by \!arclength%
  \!drawlinearsegment
  \!xS=\!xE \!yS=\!yE
  \ignorespaces}

\def\!linearsolid{%
  \!npoints=\!arclength
  \!countA=\plotsymbolspacing
  \divide\!npoints by \!countA
  \ifnum \!npoints<1 
    \!npoints=1 
  \fi
  \divide\!xdiff by \!npoints
  \divide\!ydiff by \!npoints
  \!xpos=\!xS \!ypos=\!yS
  \loop\ifnum\!npoints>-1
    \!plotifinbounds
    \advance \!xpos by \!xdiff
    \advance \!ypos by \!ydiff
    \advance \!npoints by -1
  \repeat
  \ignorespaces}

\def\!lineardashed{%
  \ifdim\!distacross>\!arclength
    \advance \!distacross by -\!arclength  
  \else
    \loop\ifdim\!distacross<\!arclength
      \!divide\!distacross\!arclength\!dimenA
      \!removept\!dimenA\!t
      \!xpos=\!t\!xdiff \advance \!xpos by \!xS
      \!ypos=\!t\!ydiff \advance \!ypos by \!yS
      \!plotifinbounds
      \advance\!distacross by \plotsymbolspacing
      \!advancedashing
    \repeat  
    \advance \!distacross by -\!arclength
  \fi
  \ignorespaces}

\def\!!advancedashing{%
  \advance\!downlength by -\plotsymbolspacing
  \ifdim \!downlength>\!zpt
  \else
    \advance\!distacross by \!downlength
    \!getnextvalueof\!uplength\from\!!UDlist
    \advance\!distacross by \!uplength
    \!getnextvalueof\!downlength\from\!!UDlist
  \fi}

\def\inboundscheckoff{%
  \def\!plotifinbounds{\!plot(\!xpos,\!ypos)}%
  \def\!initinboundscheck{\relax}\ignorespaces}
 
\inboundscheckoff
 
\def\!!plotifinbounds{%
  \ifdim \!xpos<\!checkleft
  \else
    \ifdim \!xpos>\!checkright
    \else
      \ifdim \!ypos<\!checkbot
      \else
         \ifdim \!ypos>\!checktop
         \else
           \!plot(\!xpos,\!ypos)
         \fi 
      \fi
    \fi
  \fi}

\def\!!initinboundscheck{%
  \!checkleft=\!arealloc     \advance\!checkleft by \!xorigin
  \!checkright=\!arearloc    \advance\!checkright by \!xorigin
  \!checkbot=\!areabloc      \advance\!checkbot by \!yorigin
  \!checktop=\!areatloc      \advance\!checktop by \!yorigin}

\def\!logten#1#2{%
  \expandafter\!!logten#1\!nil
  \!removept\!dimenF#2%
  \ignorespaces}

\def\!!logten#1#2\!nil{%
  \if -#1%
    \!dimenF=\!zpt
    \def\!next{\ignorespaces}%
  \else
    \if +#1%
      \def\!next{\!!logten#2\!nil}%
    \else
      \if .#1%
        \def\!next{\!!logten0.#2\!nil}%
      \else
        \def\!next{\!!!logten#1#2..\!nil}%
      \fi
    \fi
  \fi
  \!next}

\def\!!!logten#1#2.#3.#4\!nil{%
  \!dimenF=1pt 
  \if 0#1%
    \!!logshift#3pt 
  \else 
    \!logshift#2/
    \!dimenE=#1.#2#3pt 
  \fi 
  \ifdim \!dimenE<\!rootten
    \multiply \!dimenE 10 
    \advance  \!dimenF -1pt
  \fi
  \!dimenG=\!dimenE
    \advance\!dimenG 10pt
  \advance\!dimenE -10pt 
  \multiply\!dimenE 10 
  \!divide\!dimenE\!dimenG\!dimenE
  \!removept\!dimenE\!t
  \!dimenG=\!t\!dimenE
  \!removept\!dimenG\!tt
  \!dimenH=\!tt\!tenAe
    \divide\!dimenH 100
  \advance\!dimenH \!tenAc
  \!dimenH=\!tt\!dimenH
    \divide\!dimenH 100   
  \advance\!dimenH \!tenAa
  \!dimenH=\!t\!dimenH
    \divide\!dimenH 100 
  \advance\!dimenF \!dimenH}

\def\!logshift#1{%
  \if #1/%
    \def\!next{\ignorespaces}%
  \else
    \advance\!dimenF 1pt 
    \def\!next{\!logshift}%
  \fi 
  \!next}
 
 \def\!!logshift#1{%
   \advance\!dimenF -1pt
   \if 0#1%
     \def\!next{\!!logshift}%
   \else
     \if p#1%
       \!dimenF=1pt
       \def\!next{\!dimenE=1p}%
     \else
       \def\!next{\!dimenE=#1.}%
     \fi
   \fi
   \!next}

\def\beginpicture{%
  \setbox\!picbox=\hbox\bgroup%
  \!xleft=\maxdimen  
  \!xright=-\maxdimen
  \!ybot=\maxdimen
  \!ytop=-\maxdimen}
 
\def\endpicture{%
  \ifdim\!xleft=\maxdimen
    \!xleft=\!zpt \!xright=\!zpt \!ybot=\!zpt \!ytop=\!zpt 
  \fi
  \global\!Xleft=\!xleft \global\!Xright=\!xright
  \global\!Ybot=\!ybot \global\!Ytop=\!ytop
  \egroup%
  \ht\!picbox=\!Ytop  \dp\!picbox=-\!Ybot
  \ifdim\!Ybot>\!zpt
  \else 
    \ifdim\!Ytop<\!zpt
      \!Ybot=\!Ytop
    \else
      \!Ybot=\!zpt
    \fi
  \fi
  \hbox{\kern-\!Xleft\lower\!Ybot\box\!picbox\kern\!Xright}}
 
\def\endpicturesave <#1,#2>{%
  \endpicture \global #1=\!Xleft \global #2=\!Ybot \ignorespaces}

\def\setcoordinatesystem{%
  \!ifnextchar{u}{\!getlengths }
    {\!getlengths units <\!xunit,\!yunit>}}
\def\!getlengths units <#1,#2>{%
  \!xunit=#1\relax
  \!yunit=#2\relax
  \!ifcoordmode 
    \let\!SCnext=\!SCccheckforRP
  \else
    \let\!SCnext=\!SCdcheckforRP
  \fi
  \!SCnext}
\def\!SCccheckforRP{%
  \!ifnextchar{p}{\!cgetreference }
    {\!cgetreference point at {\!xref} {\!yref} }}
\def\!cgetreference point at #1 #2 {%
  \edef\!xref{#1}\edef\!yref{#2}%
  \!xorigin=\!xref\!xunit  \!yorigin=\!yref\!yunit  
  \!initinboundscheck 
  \ignorespaces}
\def\!SCdcheckforRP{%
  \!ifnextchar{p}{\!dgetreference}%
    {\ignorespaces}}
\def\!dgetreference point at #1 #2 {%
  \!xorigin=#1\relax  \!yorigin=#2\relax
  \ignorespaces}

\long\def\put#1#2 at #3 #4 {%
  \!setputobject{#1}{#2}%
  \!xpos=\!M{#3}\!xunit  \!ypos=\!M{#4}\!yunit  
  \!rotateaboutpivot\!xpos\!ypos%
  \advance\!xpos -\!xorigin  \advance\!xpos -\!xshift
  \advance\!ypos -\!yorigin  \advance\!ypos -\!yshift
  \kern\!xpos\raise\!ypos\box\!putobject\kern-\!xpos%
  \!doaccounting\ignorespaces}
 
\long\def\multiput #1#2 at {%
  \!setputobject{#1}{#2}%
  \!ifnextchar"{\!putfromfile}{\!multiput}}
\def\!putfromfile"#1"{%
  \expandafter\!multiput \input #1 /}
\def\!multiput{%
  \futurelet\!nextchar\!!multiput}
\def\!!multiput{%
  \if *\!nextchar
    \def\!nextput{\!alsoby}%
  \else
    \if /\!nextchar
      \def\!nextput{\!finishmultiput}%
    \else
      \def\!nextput{\!alsoat}%
    \fi
  \fi
  \!nextput}
\def\!finishmultiput/{%
  \setbox\!putobject=\hbox{}%
  \ignorespaces}
 
\def\!alsoat#1 #2 {%
  \!xpos=\!M{#1}\!xunit  \!ypos=\!M{#2}\!yunit  
  \!rotateaboutpivot\!xpos\!ypos%
  \advance\!xpos -\!xorigin  \advance\!xpos -\!xshift
  \advance\!ypos -\!yorigin  \advance\!ypos -\!yshift
  \kern\!xpos\raise\!ypos\copy\!putobject\kern-\!xpos%
  \!doaccounting
  \!multiput}
 
\def\!alsoby*#1 #2 #3 {%
  \!dxpos=\!M{#2}\!xunit \!dypos=\!M{#3}\!yunit 
  \!rotateonly\!dxpos\!dypos
  \!ntemp=#1%
  \!!loop\ifnum\!ntemp>0
    \advance\!xpos by \!dxpos  \advance\!ypos by \!dypos
    \kern\!xpos\raise\!ypos\copy\!putobject\kern-\!xpos%
    \advance\!ntemp by -1
  \repeat
  \!doaccounting 
  \!multiput}
 
\def\accountingon{\def\!doaccounting{\!!doaccounting}\ignorespaces}

\accountingon
\def\!!doaccounting{%
  \!xtemp=\!xpos  
  \!ytemp=\!ypos
  \ifdim\!xtemp<\!xleft 
     \!xleft=\!xtemp 
  \fi
  \advance\!xtemp by  \!wd 
  \ifdim\!xright<\!xtemp 
    \!xright=\!xtemp
  \fi
  \advance\!ytemp by -\!dp
  \ifdim\!ytemp<\!ybot  
    \!ybot=\!ytemp
  \fi
  \advance\!ytemp by  \!dp
  \advance\!ytemp by  \!ht 
  \ifdim\!ytemp>\!ytop  
    \!ytop=\!ytemp  
  \fi}
 
\long\def\!setputobject#1#2{%
  \setbox\!putobject=\hbox{#1}%
  \!ht=\ht\!putobject  \!dp=\dp\!putobject  \!wd=\wd\!putobject
  \wd\!putobject=\!zpt
  \!xshift=.5\!wd   \!yshift=.5\!ht   \advance\!yshift by -.5\!dp
  \edef\!putorientation{#2}%
  \expandafter\!SPOreadA\!putorientation[]\!nil%
  \expandafter\!SPOreadB\!putorientation<\!zpt,\!zpt>\!nil\ignorespaces}
 
\def\!SPOreadA#1[#2]#3\!nil{\!etfor\!orientation:=#2\do\!SPOreviseshift}
 
\def\!SPOreadB#1<#2,#3>#4\!nil{\advance\!xshift by -#2\advance\!yshift by -#3}
 
\def\!SPOreviseshift{%
  \if l\!orientation 
    \!xshift=\!zpt
  \else 
    \if r\!orientation 
      \!xshift=\!wd
    \else 
      \if b\!orientation
        \!yshift=-\!dp
      \else 
        \if B\!orientation 
          \!yshift=\!zpt
        \else 
          \if t\!orientation 
            \!yshift=\!ht
          \fi 
        \fi
      \fi
    \fi
  \fi}

\long\def\!dimenput#1#2(#3,#4){%
  \!setputobject{#1}{#2}%
  \!xpos=#3\advance\!xpos by -\!xshift
  \!ypos=#4\advance\!ypos by -\!yshift
  \kern\!xpos\raise\!ypos\box\!putobject\kern-\!xpos%
  \!doaccounting\ignorespaces}

\def\!setdimenmode{%
  \let\!M=\!M!!\ignorespaces}
\def\!setcoordmode{%
  \let\!M=\!M!\ignorespaces}
\def\!ifcoordmode{%
  \ifx \!M \!M!}
\def\!ifdimenmode{%
  \ifx \!M \!M!!}
\def\!M!#1#2{#1#2} 
\def\!M!!#1#2{#1}
\!setcoordmode
\let\setdimensionmode=\!setdimenmode
\let\setcoordinatemode=\!setcoordmode

\def\!stack[#1]{%
  \let\!lglue=\hfill \let\!rglue=\hfill
  \expandafter\let\csname !#1glue\endcsname=\relax
  \!ifnextchar<{\!!stack}{\!!stack<\stackleading>}}
\def\!!stack<#1>#2{%
  \vbox{\def\!valueslist{}\!ecfor\!value:=#2\do{%
    \expandafter\!rightappend\!value\withCS{\\}\to\!valueslist}%
    \!lop\!valueslist\to\!value
    \let\\=\cr\lineskiplimit=\maxdimen\lineskip=#1%
    \baselineskip=-1000pt\halign{\!lglue##\!rglue\cr \!value\!valueslist\cr}}%
  \ignorespaces}

\def\!lines[#1]#2{%
  \let\!lglue=\hfill \let\!rglue=\hfill
  \expandafter\let\csname !#1glue\endcsname=\relax
  \vbox{\halign{\!lglue##\!rglue\cr #2\crcr}}%
  \ignorespaces}

\def\!Lines[#1]#2{%
  \let\!lglue=\hfill \let\!rglue=\hfill
  \expandafter\let\csname !#1glue\endcsname=\relax
  \vtop{\halign{\!lglue##\!rglue\cr #2\crcr}}%
  \ignorespaces}

\def\setplotsymbol(#1#2){%
  \!setputobject{#1}{#2}
  \setbox\!plotsymbol=\box\!putobject%
  \!plotsymbolxshift=\!xshift 
  \!plotsymbolyshift=\!yshift 
  \ignorespaces}
 
\setplotsymbol({\fiverm .})

\def\!!plot(#1,#2){%
  \!dimenA=-\!plotxorigin \advance \!dimenA by #1
  \!dimenB=-\!plotyorigin \advance \!dimenB by #2
  \kern\!dimenA\raise\!dimenB\copy\!plotsymbol\kern-\!dimenA%
  \ignorespaces}
 
\def\!!!plot(#1,#2){%
  \!dimenA=-\!plotxorigin \advance \!dimenA by #1

  \!dimenB=-\!plotyorigin \advance \!dimenB by #2
  \kern\!dimenA\raise\!dimenB\copy\!plotsymbol\kern-\!dimenA%
  \!countE=\!dimenA
  \!countF=\!dimenB
  \immediate\write\!replotfile{\the\!countE,\the\!countF.}%
  \ignorespaces}

\def\savelinesandcurves on "#1" {%
  \immediate\closeout\!replotfile
  \immediate\openout\!replotfile=#1%
  \let\!plot=\!!!plot}

\def\dontsavelinesandcurves {%
  \let\!plot=\!!plot}
\dontsavelinesandcurves

{\catcode`\%=11\xdef\!Commentsignal{
\def\writesavefile#1 {%
  \immediate\write\!replotfile{\!Commentsignal #1}%
  \ignorespaces}

\def\replot"#1" {%
  \expandafter\!replot\input #1 /}
\def\!replot#1,#2. {%
  \!dimenA=#1sp
  \kern\!dimenA\raise#2sp\copy\!plotsymbol\kern-\!dimenA
  \futurelet\!nextchar\!!replot}
\def\!!replot{%
  \if /\!nextchar 
    \def\!next{\!finish}%
  \else
    \def\!next{\!replot}%
  \fi
  \!next}

\def\!Pythag#1#2#3{%
  \!dimenE=#1\relax                                     
  \ifdim\!dimenE<\!zpt 
    \!dimenE=-\!dimenE 
  \fi
  \!dimenF=#2\relax
  \ifdim\!dimenF<\!zpt 
    \!dimenF=-\!dimenF 
  \fi
  \advance \!dimenF by \!dimenE
  \ifdim\!dimenF=\!zpt 
    \!dimenG=\!zpt
  \else 
    \!divide{8\!dimenE}\!dimenF\!dimenE
    \advance\!dimenE by -4pt
      \!dimenE=2\!dimenE
    \!removept\!dimenE\!!t
    \!dimenE=\!!t\!dimenE
    \advance\!dimenE by 64pt
    \divide \!dimenE by 2
    \!dimenH=7pt
    \!!Pythag\!!Pythag\!!Pythag
    \!removept\!dimenH\!!t
    \!dimenG=\!!t\!dimenF
    \divide\!dimenG by 8
  \fi
  #3=\!dimenG
  \ignorespaces}

\def\!!Pythag{
  \!divide\!dimenE\!dimenH\!dimenI
  \advance\!dimenH by \!dimenI
    \divide\!dimenH by 2}

\def\placehypotenuse for <#1> and <#2> in <#3> {%
  \!Pythag{#1}{#2}{#3}}

\def\!qjoin (#1,#2) (#3,#4){%
  \advance\!intervalno by 1
  \!ifcoordmode
    \edef\!xmidpt{#1}\edef\!ymidpt{#2}%
  \else
    \!dimenA=#1\relax \edef\!xmidpt{\the\!dimenA}%
    \!dimenA=#2\relax \edef\!xmidpt{\the\!dimenA}%
  \fi
  \!xM=\!M{#1}\!xunit  \!yM=\!M{#2}\!yunit   \!rotateaboutpivot\!xM\!yM
  \!xE=\!M{#3}\!xunit  \!yE=\!M{#4}\!yunit   \!rotateaboutpivot\!xE\!yE
  \!dimenA=\!xM  \advance \!dimenA by -\!xS
  \!dimenB=\!xE  \advance \!dimenB by -\!xM
  \!xB=3\!dimenA \advance \!xB by -\!dimenB
  \!xC=2\!dimenB \advance \!xC by -2\!dimenA
  \!dimenA=\!yM  \advance \!dimenA by -\!yS%
  \!dimenB=\!yE  \advance \!dimenB by -\!yM%
  \!yB=3\!dimenA \advance \!yB by -\!dimenB%
  \!yC=2\!dimenB \advance \!yC by -2\!dimenA%
  \!xprime=\!xB  \!yprime=\!yB
  \!dxprime=.5\!xC  \!dyprime=.5\!yC
  \!getf \!midarclength=\!dimenA
  \!getf \advance \!midarclength by 4\!dimenA
  \!getf \advance \!midarclength by \!dimenA
  \divide \!midarclength by 12
  \!arclength=\!dimenA
  \!getf \advance \!arclength by 4\!dimenA
  \!getf \advance \!arclength by \!dimenA
  \divide \!arclength by 12
  \advance \!arclength by \!midarclength
  \global\advance \totalarclength by \!arclength
  \ifdim\!distacross>\!arclength 
    \advance \!distacross by -\!arclength
  \else
    \!initinverseinterp
    \loop\ifdim\!distacross<\!arclength
      \!inverseinterp
      \!xpos=\!t\!xC \advance\!xpos by \!xB
        \!xpos=\!t\!xpos \advance \!xpos by \!xS
      \!ypos=\!t\!yC \advance\!ypos by \!yB
        \!ypos=\!t\!ypos \advance \!ypos by \!yS
      \!plotifinbounds
      \advance\!distacross \plotsymbolspacing
      \!advancedashing
    \repeat  
    \advance \!distacross by -\!arclength
  \fi
  \!xS=\!xE
  \!yS=\!yE
  \ignorespaces}

\def\!getf{\!Pythag\!xprime\!yprime\!dimenA%
  \advance\!xprime by \!dxprime
  \advance\!yprime by \!dyprime}

\def\!initinverseinterp{%
  \ifdim\!arclength>\!zpt
    \!divide{8\!midarclength}\!arclength\!dimenE
    \ifdim\!dimenE<\!wmin \!setinverselinear
    \else 
      \ifdim\!dimenE>\!wmax \!setinverselinear
      \else
        \def\!inverseinterp{\!inversequad}\ignorespaces
         \!removept\!dimenE\!Ew
         \!dimenF=-\!Ew\!dimenE
         \advance\!dimenF by 32pt
         \!dimenG=8pt 
         \advance\!dimenG by -\!dimenE
         \!dimenG=\!Ew\!dimenG
         \!divide\!dimenF\!dimenG\!beta
         \!gamma=1pt
         \advance \!gamma by -\!beta
      \fi
    \fi
  \fi
  \ignorespaces}

\def\!inversequad{%
  \!divide\!distacross\!arclength\!dimenG
  \!removept\!dimenG\!v
  \!dimenG=\!v\!gamma
  \advance\!dimenG by \!beta
  \!dimenG=\!v\!dimenG
  \!removept\!dimenG\!t}

\def\!setinverselinear{%
  \def\!inverseinterp{\!inverselinear}%
  \divide\!dimenE by 8 \!removept\!dimenE\!t
  \!countC=\!intervalno \multiply \!countC 2
  \!countB=\!countC     \advance \!countB -1
  \!countA=\!countB     \advance \!countA -1
  \wlog{\the\!countB th point (\!xmidpt,\!ymidpt) being plotted 
    doesn't lie in the}%
  \wlog{ middle third of the arc between the \the\!countA th 
    and \the\!countC th points:}%
  \wlog{ [arc length \the\!countA\space to \the\!countB]/[arc length 
    \the \!countA\space to \the\!countC]=\!t.}%
  \ignorespaces}
 
\def\!inverselinear{%
  \!divide\!distacross\!arclength\!dimenG
  \!removept\!dimenG\!t}

\def\startrotation{%
  \let\!rotateaboutpivot=\!!rotateaboutpivot
  \let\!rotateonly=\!!rotateonly
  \!ifnextchar{b}{\!getsincos }%
    {\!getsincos by {\!cosrotationangle} {\!sinrotationangle} }}
\def\!getsincos by #1 #2 {%
  \edef\!cosrotationangle{#1}%
  \edef\!sinrotationangle{#2}%
  \!ifcoordmode 
    \let\!ROnext=\!ccheckforpivot
  \else
    \let\!ROnext=\!dcheckforpivot
  \fi
  \!ROnext}
\def\!ccheckforpivot{%
  \!ifnextchar{a}{\!cgetpivot}%
    {\!cgetpivot about {\!xpivotcoord} {\!ypivotcoord} }}
\def\!cgetpivot about #1 #2 {%
  \edef\!xpivotcoord{#1}%
  \edef\!ypivotcoord{#2}%
  \!xpivot=#1\!xunit  \!ypivot=#2\!yunit
  \ignorespaces}
\def\!dcheckforpivot{%
  \!ifnextchar{a}{\!dgetpivot}{\ignorespaces}}
\def\!dgetpivot about #1 #2 {%
  \!xpivot=#1\relax  \!ypivot=#2\relax
  \ignorespaces}

\def\stoprotation{%
  \let\!rotateaboutpivot=\!!!rotateaboutpivot
  \let\!rotateonly=\!!!rotateonly
  \ignorespaces}
 
\def\!!rotateaboutpivot#1#2{%
  \!dimenA=#1\relax  \advance\!dimenA -\!xpivot
  \!dimenB=#2\relax  \advance\!dimenB -\!ypivot
  \!dimenC=\!cosrotationangle\!dimenA
    \advance \!dimenC -\!sinrotationangle\!dimenB
  \!dimenD=\!cosrotationangle\!dimenB
    \advance \!dimenD  \!sinrotationangle\!dimenA
  \advance\!dimenC \!xpivot  \advance\!dimenD \!ypivot
  #1=\!dimenC  #2=\!dimenD
  \ignorespaces}

\def\!!rotateonly#1#2{%
  \!dimenA=#1\relax  \!dimenB=#2\relax 
  \!dimenC=\!cosrotationangle\!dimenA
    \advance \!dimenC -\!rotsign\!sinrotationangle\!dimenB
  \!dimenD=\!cosrotationangle\!dimenB
    \advance \!dimenD  \!rotsign\!sinrotationangle\!dimenA
  #1=\!dimenC  #2=\!dimenD
  \ignorespaces}
\def\!rotsign{}
\def\!!!rotateaboutpivot#1#2{\relax}
\def\!!!rotateonly#1#2{\relax}
\stoprotation

\def\!reverserotateonly#1#2{%
  \def\!rotsign{-}%
  \!rotateonly{#1}{#2}%
  \def\!rotsign{}%
  \ignorespaces}

\def\!getspan span <#1>{%
  \!dshade=#1\relax
  \!ifcoordmode 
    \let\!GRnext=\!GRccheckforAP
  \else
    \let\!GRnext=\!GRdcheckforAP
  \fi
  \!GRnext}
\def\!GRccheckforAP{%
  \!ifnextchar{p}{\!cgetanchor }
    {\!cgetanchor point at {\!xshadesave} {\!yshadesave} }}
\def\!cgetanchor point at #1 #2 {%
  \edef\!xshadesave{#1}\edef\!yshadesave{#2}%
  \!xshade=\!xshadesave\!xunit  \!yshade=\!yshadesave\!yunit
  \ignorespaces}
\def\!GRdcheckforAP{%
  \!ifnextchar{p}{\!dgetanchor}%
    {\ignorespaces}}
\def\!dgetanchor point at #1 #2 {%
  \!xshade=#1\relax  \!yshade=#2\relax
  \ignorespaces}

\def\setshadesymbol{%
  \!ifnextchar<{\!setshadesymbol}{\!setshadesymbol<,,,> }}

\def\!setshadesymbol <#1,#2,#3,#4> (#5#6){%
  \!setputobject{#5}{#6}%
  \setbox\!shadesymbol=\box\!putobject%
  \!shadesymbolxshift=\!xshift \!shadesymbolyshift=\!yshift
  \!dimenA=\!xshift \advance\!dimenA \!smidge
  \!override\!dimenA{#1}\!lshrinkage%
  \!dimenA=\!wd \advance \!dimenA -\!xshift
    \advance\!dimenA \!smidge
    \!override\!dimenA{#2}\!rshrinkage
  \!dimenA=\!dp \advance \!dimenA \!yshift
    \advance\!dimenA \!smidge
    \!override\!dimenA{#3}\!bshrinkage
  \!dimenA=\!ht \advance \!dimenA -\!yshift
    \advance\!dimenA \!smidge
    \!override\!dimenA{#4}\!tshrinkage
  \ignorespaces}
\def\!smidge{-.2pt}%

\def\!override#1#2#3{%
  \edef\!!override{#2}%
  \ifx \!!override\empty
    #3=#1\relax
  \else
    \if z\!!override
      #3=\!zpt
    \else
      \ifx \!!override\!blankz
        #3=\!zpt
      \else
        #3=#2\relax
      \fi
    \fi
  \fi
  \ignorespaces}
\def\!blankz{ z}

\setshadesymbol ({\fiverm .})

\def\!startvshade#1(#2,#3,#4){%
  \let\!!xunit=\!xunit%
  \let\!!yunit=\!yunit%
  \let\!!xshade=\!xshade%
  \let\!!yshade=\!yshade%
  \def\!getshrinkages{\!vgetshrinkages}%
  \let\!setshadelocation=\!vsetshadelocation%
  \!xS=\!M{#2}\!!xunit
  \!ybS=\!M{#3}\!!yunit
  \!ytS=\!M{#4}\!!yunit
  \!shadexorigin=\!xorigin  \advance \!shadexorigin \!shadesymbolxshift
  \!shadeyorigin=\!yorigin  \advance \!shadeyorigin \!shadesymbolyshift
  \ignorespaces}
 
\def\!starthshade#1(#2,#3,#4){%
  \let\!!xunit=\!yunit%
  \let\!!yunit=\!xunit%
  \let\!!xshade=\!yshade%
  \let\!!yshade=\!xshade%
  \def\!getshrinkages{\!hgetshrinkages}%
  \let\!setshadelocation=\!hsetshadelocation%
  \!xS=\!M{#2}\!!xunit
  \!ybS=\!M{#3}\!!yunit
  \!ytS=\!M{#4}\!!yunit
  \!shadexorigin=\!xorigin  \advance \!shadexorigin \!shadesymbolxshift
  \!shadeyorigin=\!yorigin  \advance \!shadeyorigin \!shadesymbolyshift
  \ignorespaces}

\def\!lattice#1#2#3#4#5{%
  \!dimenA=#1
  \!dimenB=#2
  \!countB=\!dimenB
  \!dimenC=#3
  \advance\!dimenC -\!dimenA
  \!countA=\!dimenC
  \divide\!countA \!countB
  \ifdim\!dimenC>\!zpt
    \!dimenD=\!countA\!dimenB
    \ifdim\!dimenD<\!dimenC
      \advance\!countA 1 
    \fi
  \fi
  \!dimenC=\!countA\!dimenB
    \advance\!dimenC \!dimenA
  #4=\!countA
  #5=\!dimenC
  \ignorespaces}

\def\!qshade#1(#2,#3,#4)#5(#6,#7,#8){%
  \!xM=\!M{#2}\!!xunit
  \!ybM=\!M{#3}\!!yunit
  \!ytM=\!M{#4}\!!yunit
  \!xE=\!M{#6}\!!xunit
  \!ybE=\!M{#7}\!!yunit
  \!ytE=\!M{#8}\!!yunit
  \!getcoeffs\!xS\!ybS\!xM\!ybM\!xE\!ybE\!ybB\!ybC
  \!getcoeffs\!xS\!ytS\!xM\!ytM\!xE\!ytE\!ytB\!ytC
  \def\!getylimits{\!qgetylimits}%
  \!shade{#1}\ignorespaces}
 
\def\!lshade#1(#2,#3,#4){%
  \!xE=\!M{#2}\!!xunit
  \!ybE=\!M{#3}\!!yunit
  \!ytE=\!M{#4}\!!yunit
  \!dimenE=\!xE  \advance \!dimenE -\!xS
  \!dimenC=\!ytE \advance \!dimenC -\!ytS
  \!divide\!dimenC\!dimenE\!ytB
  \!dimenC=\!ybE \advance \!dimenC -\!ybS
  \!divide\!dimenC\!dimenE\!ybB
  \def\!getylimits{\!lgetylimits}%
  \!shade{#1}\ignorespaces}
 
\def\!getcoeffs#1#2#3#4#5#6#7#8{%
  \!dimenC=#4\advance \!dimenC -#2
  \!dimenE=#3\advance \!dimenE -#1
  \!divide\!dimenC\!dimenE\!dimenF
  \!dimenC=#6\advance \!dimenC -#4
  \!dimenH=#5\advance \!dimenH -#3
  \!divide\!dimenC\!dimenH\!dimenG
  \advance\!dimenG -\!dimenF
  \advance \!dimenH \!dimenE
  \!divide\!dimenG\!dimenH#8
  \!removept#8\!t
  #7=-\!t\!dimenE
  \advance #7\!dimenF
  \ignorespaces}

\def\!shade#1{%
  \!getshrinkages#1<,,,>\!nil
  \advance \!dimenE \!xS
  \!lattice\!!xshade\!dshade\!dimenE
    \!parity\!xpos
  \!dimenF=-\!dimenF
    \advance\!dimenF \!xE
  \!loop\!not{\ifdim\!xpos>\!dimenF}
    \!shadecolumn%
    \advance\!xpos \!dshade
    \advance\!parity 1
  \repeat
  \!xS=\!xE
  \!ybS=\!ybE
  \!ytS=\!ytE
  \ignorespaces}

\def\!vgetshrinkages#1<#2,#3,#4,#5>#6\!nil{%
  \!override\!lshrinkage{#2}\!dimenE
  \!override\!rshrinkage{#3}\!dimenF
  \!override\!bshrinkage{#4}\!dimenG
  \!override\!tshrinkage{#5}\!dimenH
  \ignorespaces}
\def\!hgetshrinkages#1<#2,#3,#4,#5>#6\!nil{%
  \!override\!lshrinkage{#2}\!dimenG
  \!override\!rshrinkage{#3}\!dimenH
  \!override\!bshrinkage{#4}\!dimenE
  \!override\!tshrinkage{#5}\!dimenF
  \ignorespaces}

\def\!shadecolumn{%
  \!dxpos=\!xpos
  \advance\!dxpos -\!xS
  \!removept\!dxpos\!dx
  \!getylimits
  \advance\!ytpos -\!dimenH
  \advance\!ybpos \!dimenG
  \!yloc=\!!yshade
  \ifodd\!parity 
     \advance\!yloc \!dshade
  \fi
  \!lattice\!yloc{2\!dshade}\!ybpos%
    \!countA\!ypos
  \!dimenA=-\!shadexorigin \advance \!dimenA \!xpos
  \loop\!not{\ifdim\!ypos>\!ytpos}
    \!setshadelocation
    \!rotateaboutpivot\!xloc\!yloc%
    \!dimenA=-\!shadexorigin \advance \!dimenA \!xloc
    \!dimenB=-\!shadeyorigin \advance \!dimenB \!yloc
    \kern\!dimenA \raise\!dimenB\copy\!shadesymbol \kern-\!dimenA
    \advance\!ypos 2\!dshade
  \repeat
  \ignorespaces}
 
\def\!qgetylimits{%
  \!dimenA=\!dx\!ytC              
  \advance\!dimenA \!ytB
  \!ytpos=\!dx\!dimenA
  \advance\!ytpos \!ytS
  \!dimenA=\!dx\!ybC              
  \advance\!dimenA \!ybB
  \!ybpos=\!dx\!dimenA
  \advance\!ybpos \!ybS}
 
\def\!lgetylimits{%
  \!ytpos=\!dx\!ytB
  \advance\!ytpos \!ytS
  \!ybpos=\!dx\!ybB
  \advance\!ybpos \!ybS}
 
\def\!vsetshadelocation{
  \!xloc=\!xpos
  \!yloc=\!ypos}
\def\!hsetshadelocation{
  \!xloc=\!ypos
  \!yloc=\!xpos}

\def\!axisticks {%
  \def\!nextkeyword##1 {%
    \expandafter\ifx\csname !ticks##1\endcsname \relax
      \def\!next{\!fixkeyword{##1}}%
    \else
      \def\!next{\csname !ticks##1\endcsname}%
    \fi
    \!next}%
  \!axissetup
    \def\!axissetup{\relax}%
  \edef\!ticksinoutsign{\!ticksinoutSign}%
  \!ticklength=\longticklength
  \!tickwidth=\linethickness
  \!gridlinestatus
  \!setticktransform
  \!maketick
  \!tickcase=0
  \def\!LTlist{}%
  \!nextkeyword}

\def\ticksout{%
  \def\!ticksinoutSign{+}}

\ticksout

\def\nogridlines{%
  \def\!gridlinestatus{\!gridlinestoofalse}}
\nogridlines

\def\loggedticks{%
  \def\!setticktransform{\let\!ticktransform=\!logten}}
\def\unloggedticks{%
  \def\!setticktransform{\let\!ticktransform=\!donothing}}
\def\!donothing#1#2{\def#2{#1}}
\unloggedticks

\expandafter\def\csname !ticks/\endcsname{%
  \!not {\ifx \!LTlist\empty}
    \!placetickvalues
  \fi
  \def\!tickvalueslist{}%
  \def\!LTlist{}%
  \expandafter\csname !axis/\endcsname}

\def\!maketick{%
  \setbox\!boxA=\hbox{%
    \beginpicture
      \!setdimenmode
      \setcoordinatesystem point at {\!zpt} {\!zpt}   
      \linethickness=\!tickwidth
      \ifdim\!ticklength>\!zpt
        \putrule from {\!zpt} {\!zpt} to
          {\!ticksinoutsign\!tickxsign\!ticklength}
          {\!ticksinoutsign\!tickysign\!ticklength}
      \fi
      \if!gridlinestoo
        \putrule from {\!zpt} {\!zpt} to
          {-\!tickxsign\!xaxislength} {-\!tickysign\!yaxislength}
      \fi
    \endpicturesave <\!Xsave,\!Ysave>}%
    \wd\!boxA=\!zpt}
  
\def\!ticksin{%
  \def\!ticksinoutsign{-}%
  \!maketick
  \!nextkeyword}

\def\!ticksout{%
  \def\!ticksinoutsign{+}%
  \!maketick
  \!nextkeyword}

\def\!tickslength<#1> {%
  \!ticklength=#1\relax
  \!maketick
  \!nextkeyword}

\def\!tickslong{%
  \!tickslength<\longticklength> }

\def\!ticksshort{%
  \!tickslength<\shortticklength> }

\def\!tickswidth<#1> {%
  \!tickwidth=#1\relax
  \!maketick
  \!nextkeyword}

\def\!ticksandacross{%
  \!gridlinestootrue
  \!maketick
  \!nextkeyword}

\def\!ticksbutnotacross{%
  \!gridlinestoofalse
  \!maketick
  \!nextkeyword}

\def\!tickslogged{%
  \let\!ticktransform=\!logten
  \!nextkeyword}

\def\!ticksunlogged{%
  \let\!ticktransform=\!donothing
  \!nextkeyword}

\def\!ticksunlabeled{%
  \!tickcase=0
  \!nextkeyword}

\def\!ticksnumbered{%
  \!tickcase=1
  \!nextkeyword}

\def\!tickswithvalues#1/ {%
  \edef\!tickvalueslist{#1! /}%
  \!tickcase=2
  \!nextkeyword}

\def\!ticksquantity#1 {%
  \ifnum #1>1
    \!updatetickoffset
    \!countA=#1\relax
    \advance \!countA -1
    \!ticklocationincr=\!axisLength
      \divide \!ticklocationincr \!countA
    \!ticklocation=\!axisstart
    \loop \!not{\ifdim \!ticklocation>\!axisend}
      \!placetick\!ticklocation
      \ifcase\!tickcase
          \relax 
        \or
          \relax 
        \or
          \expandafter\!gettickvaluefrom\!tickvalueslist
          \edef\!tickfield{{\the\!ticklocation}{\!value}}%
          \expandafter\!listaddon\expandafter{\!tickfield}\!LTlist%
      \fi
      \advance \!ticklocation \!ticklocationincr
    \repeat
  \fi
  \!nextkeyword}

\def\!ticksat#1 {%
  \!updatetickoffset
  \edef\!Loc{#1}%
  \if /\!Loc
    \def\next{\!nextkeyword}%
  \else
    \!ticksincommon
    \def\next{\!ticksat}%
  \fi
  \next}    
      
\def\!ticksfrom#1 to #2 by #3 {%
  \!updatetickoffset
  \edef\!arg{#3}%
  \expandafter\!separate\!arg\!nil
  \!scalefactor=1
  \expandafter\!countfigures\!arg/
  \edef\!arg{#1}%
  \!scaleup\!arg by\!scalefactor to\!countE
  \edef\!arg{#2}%
  \!scaleup\!arg by\!scalefactor to\!countF
  \edef\!arg{#3}%
  \!scaleup\!arg by\!scalefactor to\!countG
  \loop \!not{\ifnum\!countE>\!countF}
    \ifnum\!scalefactor=1
      \edef\!Loc{\the\!countE}%
    \else
      \!scaledown\!countE by\!scalefactor to\!Loc
    \fi
    \!ticksincommon
    \advance \!countE \!countG
  \repeat
  \!nextkeyword}

\def\!updatetickoffset{%
  \!dimenA=\!ticksinoutsign\!ticklength
  \ifdim \!dimenA>\!offset
    \!offset=\!dimenA
  \fi}

\def\!placetick#1{%
  \if!xswitch
    \!xpos=#1\relax
    \!ypos=\!axisylevel
  \else
    \!xpos=\!axisxlevel
    \!ypos=#1\relax
  \fi
  \advance\!xpos \!Xsave
  \advance\!ypos \!Ysave
  \kern\!xpos\raise\!ypos\copy\!boxA\kern-\!xpos
  \ignorespaces}

\def\!gettickvaluefrom#1 #2 /{%
  \edef\!value{#1}%
  \edef\!tickvalueslist{#2 /}%
  \ifx \!tickvalueslist\!endtickvaluelist
    \!tickcase=0
  \fi}
\def\!endtickvaluelist{! /}

\def\!ticksincommon{%
  \!ticktransform\!Loc\!t
  \!ticklocation=\!t\!!unit
  \advance\!ticklocation -\!!origin
  \!placetick\!ticklocation
  \ifcase\!tickcase
    \relax 
  \or 
    \ifdim\!ticklocation<-\!!origin
      \edef\!Loc{$\!Loc$}%
    \fi
    \edef\!tickfield{{\the\!ticklocation}{\!Loc}}%
    \expandafter\!listaddon\expandafter{\!tickfield}\!LTlist%
  \or 
    \expandafter\!gettickvaluefrom\!tickvalueslist
    \edef\!tickfield{{\the\!ticklocation}{\!value}}%
    \expandafter\!listaddon\expandafter{\!tickfield}\!LTlist%
  \fi}

\def\!separate#1\!nil{%
  \!ifnextchar{-}{\!!separate}{\!!!separate}#1\!nil}
\def\!!separate-#1\!nil{%
  \def\!sign{-}%
  \!!!!separate#1..\!nil}
\def\!!!separate#1\!nil{%
  \def\!sign{+}%
  \!!!!separate#1..\!nil}
\def\!!!!separate#1.#2.#3\!nil{%
  \def\!arg{#1}%
  \ifx\!arg\!empty
    \!countA=0
  \else
    \!countA=\!arg

  \fi
  \def\!arg{#2}%
  \ifx\!arg\!empty
    \!countB=0
  \else
    \!countB=\!arg
  \fi}
 
\def\!countfigures#1{%
  \if #1/%
    \def\!next{\ignorespaces}%
  \else
    \multiply\!scalefactor 10
    \def\!next{\!countfigures}%
  \fi
  \!next}

\def\!scaleup#1by#2to#3{%
  \expandafter\!separate#1\!nil
  \multiply\!countA #2\relax
  \advance\!countA \!countB
  \if -\!sign
    \!countA=-\!countA
  \fi
  #3=\!countA
  \ignorespaces}

\def\!scaledown#1by#2to#3{%
  \!countA=#1\relax
  \ifnum \!countA<0 
    \def\!sign{-}
    \!countA=-\!countA
  \else
    \def\!sign{}%
  \fi
  \!countB=\!countA
  \divide\!countB #2\relax
  \!countC=\!countB
    \multiply\!countC #2\relax
  \advance \!countA -\!countC
  \edef#3{\!sign\the\!countB.}
  \!countC=\!countA 
  \ifnum\!countC=0 
    \!countC=1
  \fi
  \multiply\!countC 10
  \!loop \ifnum #2>\!countC
    \edef#3{#3\!zero}%
    \multiply\!countC 10
  \repeat
  \edef#3{#3\the\!countA}
  \ignorespaces}

\def\!placetickvalues{%
  \advance\!offset \tickstovaluesleading
  \if!xswitch
    \setbox\!boxA=\hbox{%
      \def\\##1##2{%
        \!dimenput {##2} [B] (##1,\!axisylevel)}%
      \beginpicture 
        \!LTlist
      \endpicturesave <\!Xsave,\!Ysave>}%
    \!dimenA=\!axisylevel
      \advance\!dimenA -\!Ysave
      \advance\!dimenA \!tickysign\!offset
      \if -\!tickysign
        \advance\!dimenA -\ht\!boxA
      \else
        \advance\!dimenA  \dp\!boxA
      \fi
    \advance\!offset \ht\!boxA 
      \advance\!offset \dp\!boxA
    \!dimenput {\box\!boxA} [Bl] <\!Xsave,\!Ysave> (\!zpt,\!dimenA)
  \else
    \setbox\!boxA=\hbox{%
      \def\\##1##2{%
        \!dimenput {##2} [r] (\!axisxlevel,##1)}%
      \beginpicture 
        \!LTlist
      \endpicturesave <\!Xsave,\!Ysave>}%
    \!dimenA=\!axisxlevel
      \advance\!dimenA -\!Xsave
      \advance\!dimenA \!tickxsign\!offset
      \if -\!tickxsign
        \advance\!dimenA -\wd\!boxA
      \fi
    \advance\!offset \wd\!boxA
    \!dimenput {\box\!boxA} [Bl] <\!Xsave,\!Ysave> (\!dimenA,\!zpt)
  \fi}

\normalgraphs
\catcode`!=12 

\catcode`@=11 \catcode`!=11
  
\let\!pictexendpicture=\endpicture 
\let\!pictexframe=\frame
\let\!pictexlinethickness=\linethickness
\let\!pictexmultiput=\multiput
\let\!pictexput=\put

\def\beginpicture{%
  \setbox\!picbox=\hbox\bgroup%
  \let\endpicture=\!pictexendpicture
  \let\frame=\!pictexframe
  \let\linethickness=\!pictexlinethickness
  \let\multiput=\!pictexmultiput
  \let\put=\!pictexput
  \let\input=\@@input   
  \!xleft=\maxdimen  
  \!xright=-\maxdimen
  \!ybot=\maxdimen
  \!ytop=-\maxdimen}

\let\frame=\!latexframe

\let\pictexframe=\!pictexframe

\let\linethickness=\!latexlinethickness
\let\pictexlinethickness=\!pictexlinethickness

\let\\=\@normalcr
\catcode`@=12 \catcode`!=12

\begingroup\makeatletter
\def\x#1#2#3#4#5#6#7\relax{\def\x{#1#2#3#4#5#6}}
\expandafter\x\fmtname xxxxxx\relax \def\y{splain}
\ifx\x\y   
\gdef\SetFigFont#1#2#3{%
  \ifnum #1<17 \tiny\else \ifnum #1<20 \small\else
  \ifnum #1<24 \normalsize\else \ifnum #1<29 \large\else
  \ifnum #1<34 \Large\else \ifnum #1<41 \LARGE\else
     \huge\fi\fi\fi\fi\fi\fi
  \csname #3\endcsname}
\else
\gdef\SetFigFont#1#2#3{\begingroup
  \count@#1\relax \ifnum 25<\count@ \count@25 \fi
  \def\x{\endgroup\@setsize\SetFigFont{#2pt}}%
  \expandafter\x
    \csname \romannumeral\the\count@ pt\expandafter\endcsname
    \csname @\romannumeral\the\count@ pt\endcsname
  \csname #3\endcsname}
\fi
\endgroup

\usepackage{amscd, amssymb, amsmath}

\newtheorem{lm}{Lemma}
\newtheorem{tm}{Theorem}

\newcommand{\proj}{\mathbf P}

\newcommand{\barr}{\overline}
\newcommand{\rarr}{\rightarrow}
\newcommand{\oh}{{\mathcal{O}}}

\newcommand{\Q}{\mathbb{Q}}
\newcommand{\Z}{\mathbb{Z}}
\newcommand{\sy}{\mathbb{S}}
\newcommand{\G}{\mathbf{G}}

\newcommand{\eqq}{\stackrel{\sim}{=}}
\newcommand{\deli}{\bigtriangleup}

\newcommand{\mm}{\overline{M}_{1,4}}

\begin{document}
\title{A Geometric Construction of Getzler's Elliptic Relation}
\author{ Rahul Pandharipande}
\date{ 12 May 1997}
\maketitle

\pagestyle{plain}
\setcounter{section}{-1}
\section{\bf{Introduction}}
Let $\mm$ be the moduli space of Deligne-Mumford
stable 4-pointed elliptic curves. 
E. Getzler has determined the natural  $\sy_4$-module structure
on the vector space $H^*(\mm, \Q)$ using modular operads and
Deligne's mixed Hodge theory [G1]. 
The dimension of the $\sy_4$-invariant space
$H^4(\mm, \Q)^{\sy_4}$ is 7. 
$\mm$ has a natural stratification by dual graph type.
Define an $\sy_4$-invariant stratum of $\mm$ to be an
$\sy_4$-orbit of closed strata of $\mm$.
The number of invariant
dimension 2 strata in $\mm$ is 9. The classes
of these invariant strata in cohomology must therefore 
satisfy at least
2 linear relations. The first relation is evident. 
Let $\deli_0$ be the boundary stratum with 
generic element corresponding to a 4-pointed nodal rational curve.
There is a natural map $\barr{M}_{0,6} \rarr \deli_0$
obtained by identifying  the markings 5 and 6.
Pushing forward the basic 
divisor linear equivalence on $\barr{M}_{0,6}$ to
$\deli_0$ yields a relation among the dimension 2 boundary
strata of $\mm$ contained in $\deli_0$. In [G2], Getzler
computes the $9 \times 9$ intersection pairing   
of the invariant strata in $\mm$. This intersection
matrix is found to have rank 7. The invariant strata therefore
span $H^4(\mm, \Q)^{\sy_4}$. The null space of the intersection
matrix is computed to find a new relation among these strata
in cohomology.

A direct construction of Getzler's new relation via 
a rational equivalence in  
the Chow group $A_2(\mm, \Q)$ is presented here. 
The idea is to equate cycles corresponding to
different degenerations of elliptic curves in the Chow
group of a space
of admissible covers. The push-forward of this equivalence to $\mm$
yields Getzler's equation.

The strategy of the construction is the following.
First,
select a point: $$\zeta=[\proj^1, p_1,p_2,p_3,p_4]\in M_{0,4}.$$
Consider the
space of degree 2 maps $\pi: E \rarr \proj^1$ 
from 4-pointed nonsingular elliptic curves to
$\proj^1$ satisfying:
\begin{enumerate}
\item[(i)] $\pi$ is not ramified at the markings of E,
\item[(ii)] The $i^{th}$ marking of $E$ lies over $p_i$.
\end{enumerate}
This space has dimension 4 and may
be compactified (after ordering the
ramification points) by the space of
pointed admissible double covers $\barr{H}_\zeta$.
These spaces are discussed in Section \ref{admi}.
There are two natural dominant maps:
$$\lambda:\barr{H}_\zeta \rarr \barr{M}_{1,4},$$
$$\rho: \barr{H}_\zeta \rarr Sym^4 (\proj^1).$$
The map $\lambda$ is the map to moduli while
$\rho([\pi])$ is defined by the degree 4 ramification divisor
of $\pi$. 
There are two natural surfaces
in $Sym^4(\proj^1)$: the loci where the generic 
divisor shapes are 2+2 and 1+3 respectively.
These surfaces are linearly equivalent (up to a scalar)
in $A_2(Sym^4(\proj^1), \Q)$ since $A_2(Sym^4 (\proj^1),\Q)=\Q$. 
The intersection 
pull-back via $\rho$ of these surfaces yields a rational
equivalence of codimension 2 cycles in 
$\barr{H}_\zeta$. This is a {\em boundary} linear equivalence because
a double cover by an nonsingular elliptic curve
is branched over 4 distinct points. The $\lambda$ push-forward
of this relation yields Getzler's relation in $A_2(\mm, \Q)$.

C. Faber has determined the Chow ring of $\barr{M}_3$ in
[Fa]. There are 12 dimension 2 boundary strata. In [Fa], the rank
$A_2(\barr{M}_3, \Q)=7$ is computed and the
strata are shown to span. Four geometric
relations are constructed among the dimension 2 boundary
strata. Faber deduces a fifth relation via 
the ring axioms and further knowledge of $A_*(\barr{M}_3, \Q)$. 
There is a map $\tau: \mm
\rarr \barr{M}_3$ obtained by identifying the marked point
1 with 2 and 3 with 4. 
The $\tau$ push-forward of Getzler's relation is easily seen
to yield a relation among the dimension 2 strata of
$\barr{M}_3$. This provided a geometric construction of
Faber's fifth relation.

Getzler's new relation yields a
differential equation for elliptic Gromov-Witten 
invariants analogous to the WDVV-equations in genus 0
(see [G2]).
Eguchi, Hori, and Xiong and, independently,
S. Katz have proposed the full (gravitational)
potential function is annihilated by a 
certain representation of the Virasoro algebra.
The Virasoro conjecture (together
with the topological recursion relations in genus 1)
predicts a remarkably simple
recursion for the elliptic Gromov-Witten invariants
of $\proj^2$ [EHX]. In the last section,
this elliptic recursion is proven via Getzler's
equation and the WDVV-equations.

The author wishes to thank D. Abramovich, P. Belorousski, C. Faber, 
E. Getzler, 
T. Graber, S. Katz, and A. Kresch for conversations about $\mm$ and related
issues.
This research was completed at the Mittag-Leffler Institute.
The author was partially supported by an NSF post-doctoral
fellowship.

\section{\bf Moduli of maps}
\label{admi}
An {\em admissible cover of degree $2$ with $4$ marked points 
and $4$ branch points} consists
of a morphism $\pi: E\rarr D$ of pointed curves
$$[E, \ P_1, \ldots,P_4],  \ \ 
[D,\ p_1, \ldots, p_4, q_a, \ldots ,q_d]$$ satisfying the
following conditions:
\begin{enumerate}
\item[(i)] $E$ is a connected, reduced, nodal curve of arithmetic genus 1. 
\item[(ii)] The markings $P_i$ lie in the nonsingular locus $E_{ns}$.
\item[(iii)] $\pi(P_i)=p_i$.
\item[(iv)] $[D, p_1, \ldots, p_4, q_a, \ldots, q_d]$ 
is an $8$-pointed  stable curve
of genus $0$.
\item[(v)] $\pi^{-1}(D_{ns})= E_{ns}$.
\item[(vi)] $\pi^{-1}(D_{sing})=E_{sing}$.
\item[(vii)] $\pi|_{E_{ns}}$ is \'etale except over the points $q_j$ where
$\pi$ is simply ramified.
\item[(viii)] If $x\in E_{sing}$, then:
\begin{enumerate}
\item[(a)] $x \in E_1 \cap E_2 $
where $E_1, E_2$ are distinct components of $E$.
\item[(b)]
$\pi(E_1), \pi(E_2)$ are distinct components of $D$.
\item[(c)]
 The ramification numbers at
$x$ of the two morphisms:
$$\pi: E_1 \rarr \pi(E_1), \ \ \ \pi:E_2 \rarr \pi(E_2)$$
are equal.
\end{enumerate}
\end{enumerate} 
Let $\barr{H}$ be the space of
4-pointed admissible degree 2 covers of $\proj^1$ 
branched at 4 points.

$\barr{H}$ is an irreducible variety
parametrizing admissible covers.
There are natural maps
\begin{equation}
\label{mmm}
\lambda:\barr{H} \rarr \barr{M}_{1,4},
\ \ \pi:\barr{H} \rarr \barr{M}_{0,8}
\end{equation}
obtained from the domain and range of the
admissible cover respectively.
The projection $\pi$ is a finite map to $\barr{M}_{0,8}$.
For each marking $i \in\{1,2,3,4\}$, there is a
natural $\Z/2\Z$-action given by switching
the sheet of the $i^{th}$-marking of $E$.
These actions induce a product action in which the
diagonal $\deli$ acts trivially. Define the group $\G$ by:
 $$\G= (\Z/2\Z)^4/\deli.$$
The action of $\G$ on $\barr{H}$ is generically
free and commutes with the projection $\pi$. 
Therefore, the quotient $\barr{H}/\G$ naturally maps to $\barr{M}_{0,8}$.
In fact, since $\barr{H}/\G \rarr \barr{M}_{0,8}$ is finite
and birational, it is an isomorphism.

Spaces of admissible covers were defined in [HM]. The
methods there may be used to construct spaces of
pointed admissible covers. A quick route to $\barr{H}$
via Kontsevich's space of stable maps will be presented here.

Let $\mu:U\rarr \barr{M}_{0,8}$ be the 8-pointed
universal curve. Let
$\barr{M}_{1,4}(\mu,2)$ be the relative space of
stable maps of double covers of the universal curve.
$\barr{M}_{1,4}(\mu,2)$ is a projective variety (see [K], [KM], [FP], [BM]).
Let $H$
denote the quasi-projective subvariety of $\barr{M}_{1,4}(\mu,2)$
corresponding to 4-pointed maps
$\mu: E \rarr D$ satisfying:
\begin{enumerate}
\item[(i)] $E$ is a nonsingular 4-pointed elliptic curve, $D$ is
           a nonsingular 8-pointed $\proj^1$. 
\item[(ii)] The markings of $E$ lie over the first four
markings of $D$.
\item[(iii)] The four ramification points of $\mu$
             lie over the last four markings of $D$.
\end{enumerate}
$H$ is an \'etale  $8$-sheeted cover of $M_{0,8}$
with a free $\G$-action.
Let $\barr{H}$ be the closure of $H$ in 
$\barr{M}_{1,4}(\mu,2)$. The points of the closure
are seen to parameterize pointed admissible covers.
The basic maps (\ref{mmm}) are obtained from the
restriction of maps defined on $\barr{M}_{1,4}(\mu,2)$.
The $\G$-action  is deduced from the universal
properties of  $\barr{M}_{1,4}(\mu,2)$.
$\barr{H}$ and $\barr{M}_{0,8}$ are 5 dimensional spaces.

Let $\zeta = [\proj^1, p_1, p_2, p_3, p_4]\in M_{0,4}$.
Let $X_\zeta$ be the fiber over $\zeta$ of the natural map
$\barr{M}_{0,8} \rarr \barr{M}_{0,4}$ determined by the first four markings.
$X_{\zeta}$ is a nonsingular 4 dimensional variety.
Another way to view $X_\zeta$ is the following.
Let $\proj^1[8]$ be the Fulton-MacPherson configuration
space [FM]. Let $\gamma_1, \ldots, \gamma_4$ be the first four
projections $\gamma_i : \proj^1[8] \rarr \proj^1$.
$X_\zeta$ is naturally identified with the
intersection:
$$\gamma_1^{-1}(p_1) \cap \gamma_2^{-1}(p_2) \cap
\gamma_3^{-1}(p_3) \cap \gamma_4^{-1}(p_4) \subset \proj^1[8].$$
In particular, the
four markings corresponding to the
ramification points yield four natural projections
$\gamma_a, \ldots, \gamma_d$ from $X_\zeta$ to
$\proj^1$. 

$X_\zeta$ is a compactification of the
space of $8$ points on $\proj^1$ where the first
four have a specified cross ratio. The boundary of
$X_\zeta$ is a divisor with normal crossings.
Moreover, $X_\zeta$ has a stratification by graph type.

Let $\barr{H}_\zeta= \pi^{-1}(X_\zeta)$.
$\barr{H}_\zeta$ is a natural compactification of the
space of 4-pointed nonsingular elliptic double covers of $\proj^1$
satisfying:
\begin{enumerate}
\item[(i)] The ramification points on $E$ are ordered and
           distinct from the markings.
\item[(ii)]  The markings lie
over a fixed 4-tuple $\zeta \in M_{0,4}$.
\end{enumerate}
$\G$ acts on $\barr{H}_\zeta$ with quotient $\barr{H}_\zeta/\G=X_\zeta$.
The restriction $\lambda: \barr{H}_\zeta \rarr \mm$
is a generically finite  dominant map of 
degree  $48$.

\section{\bf Two Surfaces}
\label{2sur}
Let $S_{2,2}$ and $S_{1,3}$ be the irreducible
surfaces in $Sym^4 (\proj^1) \eqq \proj^4$ corresponding to
divisors of shape $2+2$ and $1+3$ respectively.
The degrees of these surfaces are $4$ and $6$ respectively.
Therefore,
$$\frac{[S_{2,2}]}{4} = \frac{[S_{1,3}]}{6}$$
in $A_2(Sym^4(\proj^1))$. Chow groups here
will always be taken with $\Q$-coefficients.

Index four $\proj^1$'s by the letters $\{a,b,c,d\}$. Define: 
$$\proj^1_4= \proj^1_a \times \cdots \times \proj^1_b.$$
The following convenient notation will be used for
the diagonal subvarieties of $\proj^1_4$.
Let $(ab)$ denote diagonal where the coordinates
of $a$ and $b$ coincide. Let $(ab,cd)= (ab)\cap (cd)$.
Similarly, let $(abc)$ denote the diagonal
where the  coordinates of $a$, $b$, and $c$ coincide.
Let $\proj^1_4 \rarr Sym^4(\proj^1)$ be the $\sy_4$-quotient
map. 
The inverse image of $S_{2,2}$ in $\proj^1_4$
is  $(ab,cd)\cup(ac,bd)\cup(ad,bc)$.
Similarly, the inverse image of $S_{1,3}$
is $(abc)\cup(acd)\cup(abd)\cup(bcd)$.

\begin {lm}
\label{pbck}
Let $\mathbf{A}$ be a finite group. Let $Z$ be an irreducible
algebraic variety with an $\mathbf{A}$-action and a quotient
morphism $\alpha: Z \rarr Z/\mathbf{A}$.
There exists a pull-back $\alpha^*: A_*(Z/\mathbf{A}) \rarr
A_*(Z)$ defined by:
$$\alpha^*([V])=  {|\text{\em Stab}(V)|}\cdot
{[\alpha^{-1}(V)_{red}]}$$
where $V$ is an irreducible subvariety of $Z/\mathbf{A}$,
the scheme $\alpha^{-1}(V)_{red}$ is the reduced preimage of $V$,
and $|\text{\em Stab}(V)|$ is the size of the generic
stabilizer of points over $V$.
\end{lm}

By Lemma \ref{pbck}, there is an equality in
$A_2(\proj_4^1)$:
\begin{equation}
\label{prodd}
[ab,cd] + [ac,bd] + [ad,bc]  = 
[abc] + [acd] +  [abd] +  [bcd].
\end{equation}
This equality can of course be checked directly. 
Lemma \ref{pbck} will be used in Section \ref{tacal}. See [V]
for a proof of Lemma \ref{pbck}.

\section{Refined Intersection}
\label{reff}
The notation of Sections \ref{admi} and \ref{2sur} is followed here.
The four projections 
$\gamma_a,\ldots, \gamma_d$ from $X_\zeta$ to
$\proj^1$ yield a map: 
$$\Gamma: X_\zeta \rarr \proj^1_4.$$
The pull-back of  the relation
(\ref{prodd}) by the map $\Gamma^*: A^*(\proj^1_4) \rarr
A^*(X_\zeta)$ is now analyzed.

First, the left side (\ref{prodd}) is considered.
The scheme theoretic inverse $\Gamma^{-1}(ab)$ is
a reduced union of boundary divisors of $X_\zeta$. 
Therefore, since the boundary in $X_\zeta$ has normal crossings,
$$\Gamma^{-1}(ab,cd)= \Gamma^{-1}(ab) \cap \Gamma^{-1}(cd)$$
is a reduced union of 2 and 3 dimensional strata
$X_\zeta$. In fact, in the notation of Appendices B and C ,
$$\Gamma^{-1}(ab,cd)= S_1 \cup S_2 \cup S_2 \cup  S_3 \cup T_3 \cup T_4.$$
$S_2$ occurs with 2 different round marking assignments.
While the strata $S_i$ are of expected dimension 2,
$T_3$ and $T_4$ are of excess dimension. Excess intersection
theory yields a dimension 2 cycle supported on
$T_3$ and $T_4$.
There is a bundle sequence on the divisor $T=T_3 \cup T_4$
which is exact away from the 1 dimensional intersections
of $T$ with the 2 dimensional components of $\Gamma^{-1}(ab,cd)$:
$$0\rarr \oh_T(T) \rarr \Gamma^*(N_{ab,cd}) \rarr L \rarr 0.$$
Here, $N_{ab,cd}$ is the normal bundle of $(ab,cd)$ in $\proj^1_4$,
the first map is obtained from the differential of $\Gamma$,
and $L$ is the excess bundle. The first Chern class of
$L$ is the excess cycle on $T$ (see [Fu]).
A standard calculation yields the excess cycle:  
$$2 R_1+ \frac{1}{2} R_2+ 2 R_3 + 2 R_4 + S_4 +S_7 +S_8.$$
The intersection pull-back is thus determined by:
\begin{equation*}
\Gamma^*[ab,cd]  =   2 R_1+ \frac{1}{2} R_2+ 2 R_3 + 2 R_4  
+S_1+ 2 S_2 + S_3+ S_4 +S_7 +S_8.
\end{equation*}
The above equalities are in $A_2(X_\zeta)/K$ (see Appendix B).

The scheme theoretic inverse $\Gamma^{-1}(abc)$ is
a union of boundary divisors of $X_\zeta$:
$$\Gamma^{-1}(abc)= T_1 \cup T_2 \cup T_3 \cup T_4.$$
The excess bundle $L$ is again determined by the natural 
exact sequence on $T=\Gamma^{-1}(abc)$:
$$0\rarr \oh_T(T) \rarr \Gamma^*(N_{abc}) \rarr L \rarr 0.$$
A calculation yields the excess class in $A_2(X_\zeta)/K$:
$$
\Gamma^*[abc]  =   \frac{6}{10}R_1  + \frac{1}{2} R_2
+ \frac{13}{10}R_3 
+\frac{1}{2} R_4 
+ \frac{6}{10} R_5  + \frac{1}{10} R_6  
+\frac{1}{6} R_7 $$
$$+ S_4 + S_5 + S_6 +S_7 +S_8 +S_9. $$

The $\Gamma$ pull-back of the left
side of (\ref{prodd}) in $A_2(X_\zeta)/K$ is simply
$3\Gamma^*[ab,cd]$. Similarly, the pull-back of the
right side is $4\Gamma^*[abc]$. Hence, the
above calculation yield an equality in $A_2(X_\zeta)/K$:
\begin{equation}
\label{mainn}
\frac{36}{10}R_1  - \frac{1}{2} R_2
+ \frac{8}{10}R_3 
+ 4 R_4 
- \frac{24}{10} R_5  - \frac{4}{10} R_6  
-\frac{4}{6} R_7 = 
\end{equation}
$$-3S_1-6S_2-3S_3+S_4+4S_5+4S_6+S_7+S_8+4S_9.$$

\section{Push-Forwards}
\label{tacal}
The application of $\lambda_* \pi^*: A_2(X_\zeta) \rarr A_2(\mm)$
to relation (\ref{mainn}) yields
a boundary relation in $A_2(\mm)$.
The calculation of $\pi^*$ is obtained by Lemma 1.
It is checked the generic stabilizers of $\G$ over the
strata are trivial except in the following cases:
\begin{eqnarray*}
|\text{Stab}(R_4)| & = & 2 \\
|\text{Stab}(R_5)| & = & 2 \\
|\text{Stab}(R_6)| & = & 2 \\
|\text{Stab}(R_7)| & = & 4 \\
|\text{Stab}(S_6)| & = & 2  \\
|\text{Stab}(S_9)| & = & 2.
\end{eqnarray*}

The elements $\lambda_* \pi^*(R_i)$ are tabulated below.
See Appendix A for the notation in $A_2(\mm)$.
These elements lie in the linear span of the 4 invariant strata
of $A_2(\mm)$ with nonsingular elliptic components.  
\begin{tabbing}
$\lambda_* \pi^*(R_*)+$ \= +=+ \= $++ \deli_{*.*}$  \= + \= 
$++ \deli_{*,*}$ \= + \=
$\deli_{*,*}$ \= + \= + \= $\deli_{*,*}$ \kill
$\lambda_* \pi^*(R_1)$ \> = \> \>  \> \> 2 \> 
$\deli_{2,4}$ \> + \> 6 \>$\deli_{3,4}$\\
$\lambda_* \pi^*(R_2)$ \> = \> $96 \deli_{2,2}$ \>+ \> $32 \deli_{2,3}$
\>  \> \> \> \>\\
$\lambda_* \pi^*(R_3)$ \> = \> \> \> \> 12 \> $\deli_{2,4}$ \> \> \> \\
$\lambda_* \pi^*(R_4)$ \> = \> \> \> $16 \deli_{2,3}$
\>  \> \> \> \>\\
$\lambda_* \pi^*(R_5)$ \> = \>  \> \> \> \> \> \> 12\> $\deli_{3,4}$\\
$\lambda_* \pi^*(R_6)$ \> = \> \> \> \> 32 \> 
$\deli_{2,4}$ \> + \> 12 \>  $\deli_{3,4}$\\
$\lambda_* \pi^*(R_7)$ \> = \> \> \> $48 \deli_{2,3}$
\>  \> \> \> \> 
\end{tabbing}

The elements $\lambda_* \pi^*(S_i)$ lie in the linear span of
the 5 invariant strata contained in the boundary
divisor $\deli_0 \subset \mm$ corresponding to a 4-pointed
nodal rational curve.

\begin{tabbing}
$\lambda_* \pi^*(S_*)+$ \= +=+ \= $++* \deli_{*.*}$  \= $++* \deli_{*,*}$ \= 
$++* \deli_{*,*}$ \= $++* \deli_{*}$ \= $++* \deli_{*}$ \kill
$\lambda_* \pi^*(S_1)$ \> = \> $\frac{2}{3} \deli_{0,2}$  
\> \> \> $+2 \deli_{a}$ \> $+\frac{16}{3} \deli_{b}$ \\
$\lambda_* \pi^*(S_4)$ \> = \> $4 \deli_{0,2}$  
\> \> \> \>  \\
$\lambda_* \pi^*(S_5)$ \> = \>   
\> \> $2\deli_{0,4}$ \> \>  \\
$\lambda_* \pi^*(S_6)$ \> = \>   
\> $2 \deli_{0,3}$  \> \> \>  \\
$\lambda_* \pi^*(S_7)$ \> = \>   
\> $2 \deli_{0,3}$  \> \> \>  \\
\end{tabbing}

\noindent The push-forwards $\lambda_*\pi^*(S_2)$, $\lambda_*\pi^*(S_3)$,
$\lambda_*\pi^*(S_8)$, and $\lambda_*\pi^*(S_9)$ all vanish.

In the above formulas, the ordinary coarse moduli fundamental classes
in $\mm$ are used on the right.
The orbifold classes differ by a factor of
two for $\deli_{2,4}$, $\deli_{3,4}$, and $\deli_{0,4}$.

These push forwards are easy to compute. A representative
example will be given.
Consider the graph $R_4$:

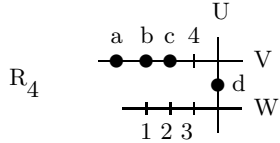
\begin{center}
  \input{r4.pictex}
\end{center}

\noindent 
with a fixed labelling of the marked points and the
round markings.
The three components of the graph have been labeled
$U$, $V$, and $W$.
The labeled graph now corresponds to an irreducible stratum $R$
of $X_\zeta$.
An admissible cover lying above this stratum consists
of the following data: 
an
elliptic double cover of $V$ ramified at $a$, $b$, $c$,
and the intersection $V \cap U$, a rational double
cover of $U$ ramified at $d$
and the intersection $V\cap W$, and an disjoint \'etale
double cover of $W$. The double cover of $W$ is a union of
two $\proj^1$'s. The markings $\{1,2,3\}$  
can be distributed in 4 ways (up to isomorphism) on the
2 components over $W$:
$$(12,3),(13,2), (23,1), (123).$$   
This yields four components of $\pi^{-1} (R)$. The component $(12,3)$
pushes forward by $\lambda$ to:

\begin{center}
  \input{ex1.pictex}
\end{center}

\noindent The $\lambda$ push-forward of the component $(123)$ lies in the
divisor stratum of $\mm$:
 
\begin{center}
  \input{ex3.pictex}
\end{center}

\noindent and is easily seen to equal:

\begin{center}
  \input{ex2.pictex}
\end{center}

\noindent Adding all the permutations and multiplying
by the automorphism factor 2 associated to $R_4$ yields
$\lambda_* \pi^*(R_4)= 16 \ \deli_{2,3}$.
\section{Getzler's relation}
The calculations of $\lambda_* \pi^*$ in Section \ref{tacal} and the
linear equivalence (\ref{mainn}) yield the relation:
\begin{equation}
\label{anss}
48 \deli_{2,2}-16 \deli_{2,3}-4 \deli_{2,4}+12\deli_{3,4}
\end{equation}
$$+2 \deli_{0,2}
+10 \deli_{0,3} +8 \deli_{0,4} -6 \deli_a -16 \deli_b=0$$
in $A_2(\mm)$. The basic linear equivalence in $A_2(\mm)$ 
obtained from $\deli_0$ is (see [G2]):
\begin{equation}
\label{ratt}
\deli_{0,2}+ 3\deli_{0,3} + 3 \deli_{0,4} -3 \deli_a -4 \deli_b=0
\end{equation}
Equation (\ref{anss}) minus twice equation (\ref{ratt}) yields:
$$4 \cdot (12 \deli_{2,2}-4 \deli_{2,3}- \deli_{2,4}+3 \deli_{3,4}
+ \deli_{0,3}+ \frac{1}{2} \deli_{0,4} - 2 \deli_b)=0$$
which is Getzler's relation [G2].
\begin{tm}
Getzler's relation in $H^4(\mm, \Q)$ is obtained from
a rational equivalence in $A_2(\mm, \Q)$.
\end{tm}

\section{Elliptic invariants of $\proj^2$}
For $d\geq 1$, let $N^{(0)}_d$ and  $N_d^{(1)}$ be the rational
and elliptic Gromov-Witten invariants of $\proj^2$.
Let $\Gamma$ and $E$ be the (quantum) rational and elliptic potentials:
\begin{eqnarray}
\Gamma(y_1,y_2) & = & 
\sum_{d\geq 1} N^{(0)}_d e^{dy_1} \frac{y_2^{3d-1}}{(3d-1)!}, \\
E(y_1,y_2) &
= & -\frac{y_1}{8} + \sum_{d\geq 3} N_d^{(1)} e^{dy_1} \frac{y_2^{3d}}
{(3d)!}.
\end{eqnarray}
The variables $y_0, y_1, y_2$ correspond to
the fundamental, line, and point classes of $\proj^2$
respectively. The elliptic sum starts in degree 3 since
$N_1^{(1)}$ and $N_2^{(1)}$ vanish.
Let
$\tilde{E}= E+ \frac{1}{8} y_1$.
The composition axiom for Gromov-Witten invariants ([KM], [RT], [BM])
and Getzler's relation immediately yield
partial differential equations satisfied by the 
elliptic potential.
Getzler uses these to obtain a complete recursion
determining the numbers $N_d^{(1)}$ [G2].
The numbers $N_d^{(1)}$ have been obtained via
more classical techniques in algebraic geometry
in [CH], [R].

The
relation predicted by Eguchi, Hori,
and Xiong from the Virasoro conjecture is:
\begin{equation}
\label{egu}
\frac{N_d^{(1)}}{(3d-1)!} = \frac{1}{12} \binom {d}{3} \frac{ N_d^{(0)} }
{(3d-1)!} +
\sum_{d_1+d_2=d, d_i>0} \frac{3d_1^2 d_2-2d_1d_2}{9}
\frac{N^{(0)}_{d_1}}{(3d_1-1)!} 
\frac{N_{d_2}^{(1)}}{(3d_2)!}.
\end{equation}
While the equation suggests divisor geometry, no
proof via such an approach is known to the author.
The relation is established here via Getzler's equation and the
WDVV-equations.

The first step is to rewrite the recursion (\ref{egu})
in differential form:
$$\tilde{E}_2 =
\frac{\Gamma_{111}-3\Gamma_{11}+2\Gamma_{1}}{72} + 
\frac{\Gamma_{11}\tilde{E}_1}{3}
-\frac{2\Gamma_1 \tilde{E}_1}{9}.$$
The subscripts denote partial differentiation by $y_1$ and
$y_2$ respectively.
Next, the equation
$$ \tilde{E}_1= \frac{y_2}{3}\tilde{E}_2 $$
is used to obtain:
$$ \tilde{E}_2 \ (1-\frac{1}{9}y_2 \Gamma_{11}+\frac{2}{27}y_2
\Gamma_1)
=\frac{\Gamma_{111}-3\Gamma_{11}+2\Gamma_{1}}{72}.$$
Finally, we see:
\begin{equation}
\label{www}
E= -\frac{y_1}{8}+ \int \frac{{\Gamma_{111}
-3\Gamma_{11}+2\Gamma_{1}}}{72}
{(1-\frac{1}{9} y_2 \Gamma_{11}+ \frac{2}{27} y_2 
\Gamma_1)^{-1}} \ dy_2.
\end{equation}
Formula (\ref{www}) is equivalent to (\ref{egu}).

In order to prove (\ref{egu}), we consider a differential
equation obtained from Getzler's relation. Let
$$\pi: \barr{M}_{1,4+3d-6}(\proj^2,d) \rarr \barr{M}_{1,4}$$
be the natural projection. Let
$$e_i: \barr{M}_{1,4 +3d-6} \rarr \proj^2$$
be the natural evaluation maps.
The codimension 2 Gromov-Witten class
$$\pi_* \big( \ [\barr{M}_{1,4+3d-6}(\proj^2,d)]^{Vir} \cap
\prod_{i=1}^{3d-2} e_i^*([\text{point}]) \ \big)$$
intersected with Getzler's relation yields the following
differential equation:
\begin{equation}
\label{fff}
36E_{11} \Gamma_{122}^2 -48 E_{12}\Gamma_{112}\Gamma_{122}
-48E_{22}\Gamma_{222}-12 E_1 \Gamma_{1122} \Gamma_{122} 
\end{equation}
$$
+24 E_1 \Gamma_{112}\Gamma_{1222} +24 E_2 \Gamma_{2222}
+2 \Gamma_{1222}\Gamma_{1112}$$ $$+ \frac{1}{2}
\Gamma_{12222}\Gamma_{111}+ \frac{3}{2} \Gamma_{22222}
-3 \Gamma_{1122}^2 =0$$
It is easily seen by examining the coefficient
recursions that (\ref{fff}) uniquely determines
the potential $E$ in form (7) from $\Gamma$.  
Therefore,  to establish (\ref{egu}), it suffices to prove 
$(\ref{www})$ determines a solution of (\ref{fff}).

Replace the terms $E_1$, $E_{11}$, $E_{12}$ 
in (\ref{fff}) by:
\begin{eqnarray*}
E_1 & = & \frac{y_2}{3} E_2 - \frac{1}{8}, \\
E_{11} & = & \frac{y_2^2}{9} E_{22} + \frac{y_2}{9} E_2, \\
E_{12} & = & \frac{y_2}{3} E_{22} + \frac{1}{3} E_2.
\end{eqnarray*}
Equation (\ref{www}) then may be used to
replace the terms $E_2$ and $E_{22}$ by
derivatives of $\Gamma$. A 
partial differential equation for $\Gamma$ is then obtained.
Equation (\ref{www}) defines a solution of (\ref{fff})  
if and only if
$\Gamma$ satisfies this differential equation.

The WDVV-equation for $\Gamma$ is:
\begin{equation}
\label{wd}
\Gamma_{222} -\Gamma^2_{112} + \Gamma_{111}\Gamma_{122}=0.
\end{equation}
$\Gamma$ also
satisfies:
\begin{equation}
\Gamma_1= \frac{y_2}{3} \Gamma_2 + \frac{1}{3} \Gamma.
\label{ll}
\end{equation}
A simple check in symbolic algebra shows the
required differential equation for $\Gamma$ is
implied by the WDVV-equation and its first two
$y_2$-derivatives together with (\ref{ll}). 
In fact, the required differential equation for $\Gamma$ 
modulo equation (\ref{ll}) is seen to be of the form:
\begin{equation}
\label{fform}
A \Psi + B \Psi_2+ C \Psi^2_{2} + D \Psi_{22} =0.
\end{equation}
$A$ , $B$, $C$, and $D$ are  
polynomials in $y_2$, $\Gamma$, and the derivatives
of $\Gamma$.  $\Psi$ is the left side of the WDVV-relation
(\ref{wd}), and the subscripts denote partial differentiation
by $y_2$.
\begin{eqnarray*}
A & = & -2^7 \cdot
(36y_2\Gamma_{22} +
 24y_2^2\Gamma_{222}- 4\Gamma_{1}\Gamma_{11}
              -9 \Gamma^2_{11}+\frac{4}{3}\Gamma_{1}^2 \\
& & 
\ \ \ \ \ \ \ \ \ +12 \Gamma_1\Gamma_{111}
-18 \Gamma_{11}\Gamma_{111} +27 \Gamma_{111}^2) \\
B & = & -2^6 \cdot (-60 \Gamma +375 \Gamma_1 -\frac{1521}{2} \Gamma_{11}
+\frac{783}{2} \Gamma_{111} \\
& & \ \ \ \ \ \ \ \ \ +2y_2 \Gamma_1^2-6 y_2\Gamma_1\Gamma_{11}+8 y_2 \Gamma_1
\Gamma_{111}  \\
& & \ \ \ \ \ \ \ \ \ -\frac{11}{2}y_2\Gamma_{11}^2-9y_2\Gamma_{11}\Gamma_{111}
+\frac{9}{2}y_2\Gamma_{111}^2) \\
C & = & -2^6 \cdot y_2^4 \\
D & = & -2^5 \cdot(y_2 \Gamma_{111}+9) \cdot
(3y_2\Gamma_{11}-2y_2\Gamma_1-27)
\end{eqnarray*}
This concludes the proof of recursion (\ref{egu}).

\begin{tm}
The Eguchi-Hori-Xiong formula holds for the elliptic Gromov-Witten
invariants of $\proj^2$.
\end{tm}

\section*{Appendix A \ \  Strata of $\mm$ in dimension 2} 
Getzler's notation for the complete set of
$\sy_4$-invariant  dimension 2 strata of $\mm$ 
is followed:

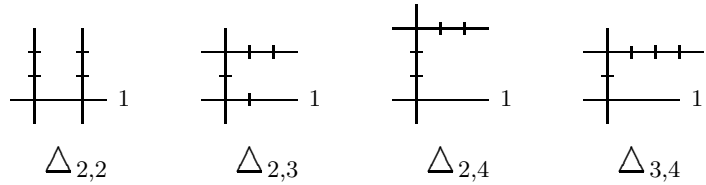
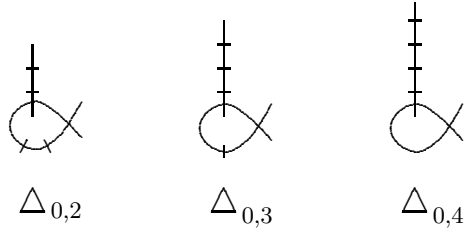
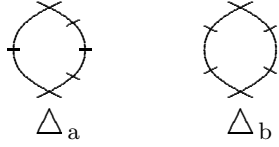
\begin{center}
  \input{g1.pictex}
\end{center}
\vspace{+5pt}
\begin{center}
  \input{g2.pictex}
\end{center}
\vspace{+5pt}
\begin{center}
\input{g3.pictex}
\end{center}

\noindent  
The short bars indicate the marked points $\{1,2,3,4\}$.
Nonsingular genus 1 components (labeled by 1) 
occur in the first row.
Each graph denotes the sum of strata in the corresponding
$\sy_4$-orbit. For example, $\deli_{2,2}$ and
$\deli_{2,3}$ are sums of 3 and 12 strata
respectively.

\section*{Appendix B \ \ Strata of $X_\zeta$ in dimension 2}
A graph $G$ below {\em together} with an
assignment $\mu$ of the round markings to the set
$\{ a,b,c,d \}$ corresponds to
an $\sy_4$-invariant dimension 3 stratum of
$\barr{M}_{0,8}$. The 8 point marking set
is $\{ 1,2,3,4,a,b,c,d \}$ as in Section \ref{admi},
and the $\sy_4$-action is on the first four markings.
By intersecting these strata of $\barr{M}_{0,8}$
with $X_\zeta$, a sum of dimension 2
strata of $X_\zeta$ is associated to $(G, \mu)$. 

For example, consider $R_3$. After an assignment
of the round markings to the set $\{ a,b,c,d\}$, there are
4 strata of $\barr{M}_{0,8}$ in the $\sy_4$-orbit. 
These yield a sum of four
strata of $X_\zeta$. Similarly, $S_3$ together
with an assignment of the round markings yields
a sum 12 strata of $X_\zeta$. 

The assignment $\mu$ of the round markings 
will be suppressed since the calculation of the map
$$\lambda_* \pi^*: A_2(X_\zeta) \rarr A_2(\mm)$$
does not depend upon $\mu$.  
Let $K=\text{Ker}(\lambda_* \pi^*)$.
Let $\mu$ and  $\mu'$ be 
two round marking assignments for a fixed graph $G$. The
two cycles corresponding to $(G, \mu)$ and $(G, \mu')$
are equal  in $A_2(X_\zeta)/K$.

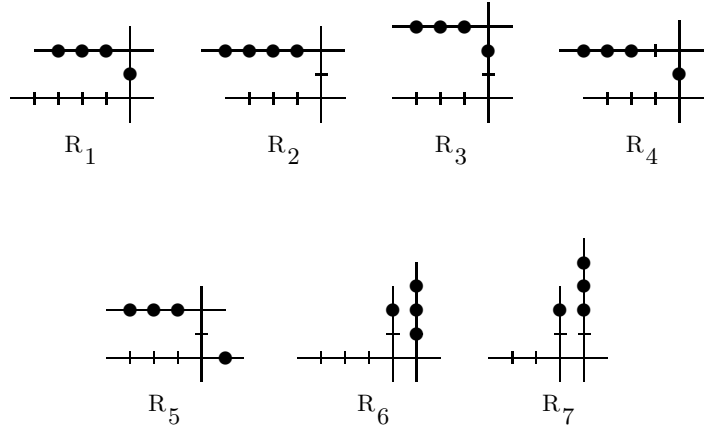
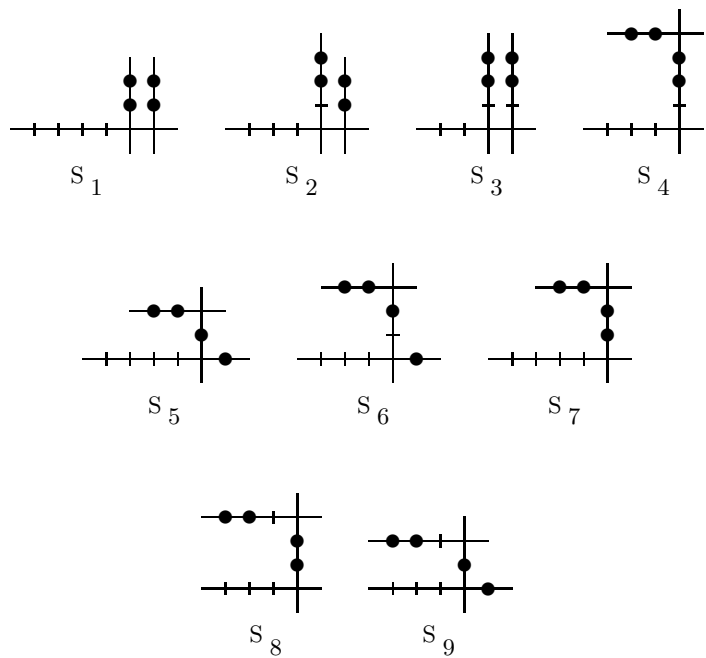
\begin{center}
\input{r1.pictex}
\end{center}
\vspace{+5pt}
\begin{center}
\input{r2.pictex}
\end{center}
\vspace{+5pt}
\begin{center}
\input{s1.pictex}
\end{center}
\vspace{+5pt}
\begin{center}
\input{s2.pictex}
\end{center}
\vspace{+5pt}
\begin{center}
\input{s3.pictex}
\end{center}

\noindent This is not a complete list. Only the strata
which are required for the computations in Section \ref{reff}
are given.

\section*{Appendix C \ \ Strata of $X_\zeta$ in dimension 3}
As in Appendix B, a
graph below together with an assignment of the
round markings to the set $\{a,b,c,d\}$ corresponds to
a sum of dimension 3 strata of $X_\zeta$.
\begin{center}
  \input{t.pictex}
\end{center}

\vspace{+10 pt}
\noindent
Department of Mathematics \\
University of Chicago \\
5734 S. University Ave. 60637 \\
rahul@math.uchicago.edu
\end{document}